\documentclass[a4paper,twoside,10pt]{article}
\usepackage[a4paper,left=3cm,right=3cm, top=3cm, bottom=3cm]{geometry}

\usepackage{amsfonts}
\usepackage{mathrsfs}
\usepackage{amsmath}
\usepackage{amssymb}
\usepackage{amsthm}
\usepackage{cite}
\usepackage{color}
\usepackage{cancel}
\usepackage{appendix}
\usepackage[cal=boondox]{mathalpha}
\usepackage{tikz}
\usepackage{soul}

\theoremstyle{plain}
\newtheorem{theorem}{Theorem}[section]
\newtheorem{corollary}[theorem]{Corollary}
\newtheorem{lemma}[theorem]{Lemma}

\theoremstyle{remark}
\newtheorem{remark}{Remark}

\newcommand{\blue}[1]{{\color{blue}{#1}}}

\newcommand{\orange}[1]{{\color{orange}{#1}}} 
\newcommand{\red}[1]{{\color{red}{#1}}} 
\newcommand{\eremk}{\hbox{}\hfill\rule{0.8ex}{0.8ex}}
\DeclareTextCompositeCommand{\u}{T1}{i}{\u\imath}

\let\div\relax
\DeclareMathOperator{\div}{\nabla\cdot}
\DeclareMathOperator{\divbf}{\nablabold\cdot}

\DeclareMathOperator{\nablabold}{\boldsymbol{\nabla}}
\newcommand{\h}{h}
\newcommand{\sigmabold}{\boldsymbol \sigma}

\newcommand{\vbftilde}{\widetilde{\vbf}}
\newcommand{\E}{K}
\newcommand{\hE}{\h_\E}

\newcommand{\F}{F}
\newcommand{\T}{T}
\newcommand{\hF}{\h_\F}

\newcommand{\FcalE}{\mathcal F^\E}

\newcommand{\Nbb}{\mathbb N}

\newcommand{\Rbb}{\mathbb R}

\newcommand{\nbf}{\mathbf n}
\newcommand{\nbfE}{\nbf_\E}
\newcommand{\nbfT}{\nbf_\T}
\newcommand{\Norm}[1]{{\left\|{#1} \right\|}}
\newcommand{\SemiNorm}[1]{{\left|{#1} \right|}}
\newcommand{\jump}[1]{\left[\!\left[#1\right]\!\right]}
\newcommand{\Normth}[1]{{\left\vert\kern-0.25ex\left\vert\kern-0.25ex\left\vert #1 
    \right\vert\kern-0.25ex\right\vert\kern-0.25ex\right\vert}}
\newcommand{\taun}{\mathcal T_n}

\newcommand{\ostar}{1^*}
\newcommand{\ubf}{\mathbf u}
\newcommand{\vbf}{\mathbf v}
\newcommand{\Vbf}{\mathbf V}

\newcommand{\Brho}{B_\rho}
\newcommand{\rhotilde}{\widetilde\rho}
\newcommand{\Brhotilde}{B_{\rhotilde}}

\newcommand{\hOmega}{\h_\Omega}
\newcommand{\hOmegatilde}{\h_{\Omegatilde}}

\newcommand{\ftilde}{\widetilde f}
\newcommand{\p}{p}
\newcommand{\q}{q}
\newcommand{\qprime}{\q'}
\newcommand{\qsharp}{\q^\sharp}
\newcommand{\qstar}{\q^*}
\newcommand{\pstar}{\p^*}
\newcommand{\psharp}{\p^\sharp}
\newcommand{\pprime}{\p'}
\newcommand{\s}{s}
\newcommand{\sstar}{\s^*}
\newcommand{\ssharp}{\s^\sharp}

\newcommand{\vh}{v_\h}
\newcommand{\Fcaln}{\mathcal F_n}
\newcommand{\FcalnI}{\mathcal F_n^I}
\newcommand{\FcalnB}{\mathcal F_n^B}
\newcommand{\FcalnD}{\mathcal F_n^D}
\newcommand{\FcalnN}{\mathcal F_n^N}

\newcommand{\nbfF}{\mathbf n_\F}
\newcommand{\Omegatilde}{\widetilde\Omega}
\newcommand{\Omegaext}{\Omega_{\text{ext}}}
\newcommand{\GammaN}{\Gamma_N}
\newcommand{\GammaD}{\Gamma_D}
\newcommand{\Gammatilde}{\widetilde\Gamma}
\newcommand{\GammatildeN}{\widetilde\Gamma_N}
\newcommand{\GammatildeD}{\widetilde\Gamma_D}
\newcommand{\xbf}{\mathbf x}
\newcommand{\xbftilde}{\widetilde{\xbf}}
\newcommand{\fbarOmega}{\overline f_\Omega}

\newcommand{\Scalstar}{\mathcal S^*}

\newcommand{\zetabold}{\boldsymbol \zeta}
\newcommand{\Psibold}{\boldsymbol \Psi}
\newcommand{\psibold}{\boldsymbol \psi}
\newcommand{\psiboldz}{\psibold^0}
\newcommand{\psiboldperp}{\psibold^\perp}
\newcommand{\Deltabold}{\boldsymbol \Delta}

\newcommand{\Wopzd}{\Wbf_0^{1,p}(\Omega)}
\newcommand{\Lp}{\L^\p(\Omega)}

\newcommand{\Lpprime}{\L^{\pprime}(\Omega)}
\newcommand{\Lpprimez}{\L^{\pprime}_0(\Omega)}
\newcommand{\nablaboldh}{\nablabold_\h}
\newcommand{\postar}{\p \ostar}
\newcommand{\Hbf}{\mathbf H}
\let\L\relax
\DeclareMathOperator{\L}{L}
\DeclareMathOperator{\W}{W}
\newcommand{\Lbf}{\mathbf L}
\newcommand{\Lbb}{\underline{\mathbf L}}
\newcommand{\Wbf}{\mathbf W}

\newcommand{\boldalpha}{\boldsymbol{\alpha}}
\newcommand{\CSob}[4]{C_{\text{Sob}}(#1,#2,#3,#4)}
\newcommand{\CTR}[4]{C_{TR}(#1,#2,#3,#4)}

\newcommand{\CG}[2]{C_{G}(#1,#2)}
\newcommand{\CPS}[2]{C_{PS}(#1,#2)}
\newcommand{\FET}{\F_\E^\T}
\newcommand{\frakTE}{\mathfrak T_\E}
\newcommand{\tildefrakTE}{\widetilde{\mathfrak T}_\E}
\newcommand{\rhoGamma}{\rho_{\Gamma}}
\newcommand{\PizFcaln}{\Pi^{0}_{\Fcaln}}
\newcommand{\Piztaun}{\Pi^{0}_{\taun}}
\DeclareMathOperator{\card}{card}
\newcommand{\htildeF}{\widetilde{\h}_\F}
\newcommand{\CBA}[3]{C_{BA}(#1,#2,#3)}

\newcommand{\CNL}[3]{C_{NL}(#1,#2,#3)}
\newcommand{\betaz}[3]{\beta(#1,#2,#3)}
\newcommand{\phinu}{\varphi_\nu}
\newcommand{\omeganu}{\omega_\nu}

\newcommand{\omegatildej}{\widetilde\omega_j}
\newcommand{\vbfnu}{\vbf_\nu}
\newcommand{\gammatilde}{\widetilde\gamma}
\newcommand{\Nfraktilde}{\widetilde{\mathfrak N}}
\newcommand{\Gcalp}{\underline{\mathbfcal{G}}_\p}
\newcommand{\Lcalp}{\mathbfcal{L}_\p}
\newcommand{\Dcalp}{\mathcal D_\p}
\newcommand{\Scalp}{\mathbfcal{S}_\p}
\newcommand{\Scalpstar}{\mathcal S_\p^*}
\newcommand{\nablau}{\underline{\nablabold}}

\newcommand{\Vbfz}{\Vbf^0}
\newcommand{\Vbfpperp}{\Vbf^\perp_\p}
\newcommand{\Vbfpz}{\Vbf^0_\p}
\newcommand{\wbf}{\mathbf w}
\let\H\relax
\DeclareMathOperator{\H}{H}

\title{\scriptsize Sobolev--Poincar\'e inequalities for piecewise $\W^{1,\p}$ functions over general polytopic meshes}
\author{\scriptsize Michele Botti\thanks{MOX, Department of Mathematics, Politecnico di Milano, 20133 Milano, Italy (michele.botti@polimi.it)},   
\and \scriptsize Lorenzo Mascotto\thanks{Department of Mathematics and Applications, University of Milano-Bicocca, 20125 Milan, Italy (lorenzo.mascotto@unimib.it); 
Faculty of Mathematics, University of Vienna, 1090 Vienna, Austria;
IMATI-CNR, 27100 Pavia, Italy}}
\date{}

\begin{document}

\maketitle

\begin{abstract}
\noindent
We establish Sobolev--Poincar\'e inequalities
for piecewise $\W^{1,\p}$ functions
over families of fairly general polytopic
(thence also shape-regular simplicial and Cartesian)
meshes in any dimension;
amongst others, they cover the case of
standard Poincar\'e inequalities for piecewise $\W^{1,\p}$ functions
and can be useful in the analysis of nonconforming finite element 
discretizations of nonlinear problems.
Crucial tools in their derivation
are novel Sobolev--trace inequalities
and $\W^{1,\p}$-stable right-inverses
of the divergence satisfying
mixed boundary conditions.
We provide estimates with constants
having an explicit dependence on
the geometric properties of the domain
and the underlying family of polytopic meshes.



\end{abstract}

\begin{flushright}
{\footnotesize\textsl{\indent Lieber ein Ende mit Schrecken
als ein Schrecken ohne Ende}}
\end{flushright}

\section{Introduction}

\paragraph*{State-of-the-art: general framework.}
Sobolev--Poincar\'e inequalities for piecewise $\W^{1,p}$ functions
are an essential tool in the analysis of several
nonconforming methods
and have been the objective of extensive studies
over the last four decades.
In the literature, the adjectives
\emph{broken}, \emph{piecewise},
and \emph{discrete} are typically associated
with this type of inequalities.
On occasions, we shall adopt either of these nomenclatures.

\paragraph*{State-of-the-art: functions in nonconforming Sobolev spaces.}
Arnold~\cite{Arnold:1982} proved
a broken Poincar\'e inequality for piecewise $\H^1$
functions in two dimensions on triangular meshes,
using elliptic regularity results in nonsmooth domains.
About 20 years later and undertaking a different avenue,
Brenner~\cite{Brenner:2003} extended Arnold's results 
to three dimensions and more general meshes;
the results in that work are based on compact embedding
arguments and the approximation properties
of certain DoFs-averaging operators;
a similar analysis was carried out for $\H^2$-type
inequalities a year later~\cite{Brenner-Want-Zhao:2004}.
Analogous results were derived
by Knobloch \cite{Knobloch:2001}
based on the stability and approximation
properties of Cl\'ement-type operators
(under more restrictive assumptions on the mesh)
and Vorhal\'ik \cite{Vohralik:2005}
only based on local Poincar\'e inequalities
and arithmetic averaging properties;
this latter reference is particularly remarkable
as it is one of the few containing
explicit dependence on several parameters
in the constant.
Lasis and S\"uli~\cite{Lasis.Suli:2003}
extended with similar arguments
the results of the papers mentioned above
to the case of Sobolev--Poincar\'e inequalities
for piecewise $\H^1$ functions
for partitions of simplices and affine maps of hypercubes.
Poincar\'e inequalities for piecewise $\W^{1,p}$ functions
in two and three dimensions
can be found in \cite[Sect. 10.6]{Brenner-Scott:2008}.

\paragraph*{State-of-the-art: piecewise polynomials.}
Broken Poincar\'e inequalities
for nonconforming methods based on polynomial spaces
were studied in even earlier works:
Temam~\cite{Temam:2001} proved
such inequalities for linear Crouzeix-Raviart spaces;
Thomas extended that analysis to
higher order Crouzeix-Raviart spaces~\cite{Thomas:1977};
see also the work by Stummel in~\cite{Stummel:1980}
dealing with other nonconforming methods.
The analysis of such references hinged upon contradiction
and finite dimensional arguments.
Similar results in the context of finite volumes
were discussed later in \cite[Section 4.3]{Coudiere-Gallouet-Herbin:2001};
in particular, the inequality therein was proven
for piecewise constant functions.
Sobolev inequalities involving discrete $\W^{1,p}$ norms
for finite volumes were detailed in
\cite[Section 5.1]{Eymard.Gallouet.ea:2010},
\cite{Glitzky-Griepentrog:2010},
and \cite[Sections~3 and~4]{Bessemoulin-Chatard.ea:2014}.
Buffa and Ortner~\cite{Buffa.Ortner:2009},
Di Pietro and Ern
\cite[Theorem 6.1]{Di-Pietro.Ern:2010},
and Droniou and Yemm \cite[Section~3]{Droniou-Yemm:2024}
proved similar inequalities for general order
piecewise polynomials over polytopic meshes
in two and three dimensions satisfying
standard \cite{Buffa.Ortner:2009, Di-Pietro.Ern:2010}
as well as more general~\cite{Droniou-Yemm:2024}
regularity assumptions;
in all three cases, the proof hinged upon the use
of polynomial trace inverse inequalities.
Related inequalities
in the framework of the staggered DG method
can be found in~\cite{Zhao-Chung-Park-Zhou:2021}.
Finally, discrete Poincar\'e inequalities
in~$\Hbf(\text{div})$ and~$\Hbf(\text{curl})$
were recently investigated in~\cite{Ern-Guzman-Potu-Vohralik:2025}.

\paragraph*{Main results and goals of the paper.}
We prove Sobolev--Poincar\'e inequalities
generalizing those in the references above
in several respects:
\begin{enumerate}
\item[(i)] the inequalities are proven
for piecewise functions over fairly general
polytopic meshes in any dimension in the sense
of, e.g., \cite[Assumption 2.1]{Cangiani-Dong-Georgoulis:2016},
see Section~\ref{subsection:meshes-broken}
below for a precise description of such meshes;
\item[(ii)] the constants appearing in the estimates
are fully explicit in terms of the geometry of the mesh
and constants appearing in direct estimates
on domains that are star-shaped with respect to a ball,
while, on domains that can be decomposed into
shape-regular simplicial meshes, the constants
have an explicit dependence
on the corresponding
shape-regularity parameter;
\item[(iii)] the employed broken Sobolev norm
contains an arbitrarily small number of boundary terms;
\item[(iv)] all the results are established in broken
$\W^{1,p}$ spaces and not for piecewise polynomials
or other finite dimensional spaces,
which is instead the avenue undertaken in
\cite{Buffa.Ortner:2009, Di-Pietro.Ern:2010, Eymard.Gallouet.ea:2010, Bessemoulin-Chatard.ea:2014, Glitzky-Griepentrog:2010};
the norms appearing in the estimates
are based on the maximal Lebesgue regularity
differently from \cite[Sect. 10.6]{Brenner-Scott:2008},
where only $\p$-type norms are used.
\end{enumerate}
Each of these results has important consequences:
\begin{enumerate}
\item[(i)] these inequalities can be used in the analysis
of several polytopic methods, in particular those based on families
of polytopic meshes with an arbitrary number of facets;
\item[(ii)] these inequalities are instrumental
in proving the well-posedness and convergence analysis
of polytopic methods, since the corresponding
involved constants are fully explicit
with respect to the geometric properties
of the family of meshes,
cf. Remark~\ref{remark:applicability} below;
\item[(iii)] these inequalities allow for the analysis of problems
with mixed boundary conditions;
\item[(iv)] they are useful for the analysis of nonconforming
methods for nonlinear problems
and Galerkin methods that are not polynomial based
(such as Trefftz and extended methods).
\end{enumerate}

\paragraph*{Important tools in the analysis.}
The proposed analysis neither relies on enriching operators,
as instead done in
\cite{Lasis.Suli:2003, Brenner:2003, Buffa.Ortner:2009},
nor on continuous embeddings of the space of functions
with bounded variation, as done in \cite{Eymard.Gallouet.ea:2010, Bessemoulin-Chatard.ea:2014,Di-Pietro.Ern:2010}.
Rather, it requires the generalization of
two fundamental results in the Sobolev spaces theory:
continuous Sobolev--trace
and  Babu\v ska--Aziz inequalities.
\medskip

We extend local continuous trace inequalities as in,
e.g., \cite[Lemma 12.15]{Ern-Guermond:2021}
and \cite[Lemma 1.31]{DiPietro-Droniou:2020},
to Sobolev--trace inequalities that are based on
more general meshes, Sobolev norms,
and maximal Lebesgue regularity;
see Section~\ref{section:1st-ingredient} below.
\medskip

Babu\v ska--Aziz inequalities express the stability of right-inverses
for several differential operators
and are strictly related to other results
such as Ne\v cas-Lions inequalities,
inf-sup conditions,
and bounds on the spectrum of the Cosserat operator;
see, e.g., \cite{Costabel-Dauge:2015} for the Hilbertian case
and Appendix~\ref{appendix:homogeneous-BCs-identities-Banach}
for the Banach case.
In the literature, the construction
of a right-inverse of the divergence
is done at least in two different ways:
one~\cite{Arnold-Scott-Vogelius:1988}
is based on solving suitable elliptic problems
and exploiting the stability
of extension operators in Sobolev spaces;
the other \cite{Bogovskii:1979, Guzman-Salgado:2021}
is based on using integral operators.
Such right-inverses are typically endowed with
homogeneous Dirichlet boundary conditions
or without boundary conditions:
the corresponding integral operators
go under the name of Bogovskiĭ and
generalized Poincar\'e operators.
However, we are not aware of
the $\W^{1,\p}$-version Babu\v ska--Aziz inequalities
for the case of boundary conditions
assigned on part of the boundary of the domain,
with constants that are explicit
in terms of the shape of the domain.
Therefore, in Section~\ref{section:2nd-ingredient} below,
we shall prove such inequalities with explicit constants
under certain regularity properties of the domain and its boundary.
A result in this direction
can be found in \cite[Lemma~$3.1$ and Appendix~A]{BeiraodaVeiga-Canuto-Nochetto-Vacca:2021},
which considers the Hilbertian case in two dimensions
without explicit constants.

\paragraph*{Outline of the remainder of the section
and list of the relevant results.}
In Section~\ref{subsection:domains-Sobolev}, we introduce
some domains of interest, their geometric properties,
and establish the notation for Sobolev spaces;
Section~\ref{subsection:meshes-broken} is devoted to introduce
families of polytopic meshes
and corresponding broken Sobolev spaces;
the two main results of the paper,
whose proofs require two technical tools discussed
in Section~\ref{subsection:main-technical-tools},
are discussed in Section~\ref{subsection:main-results}.
For the reader's convenience,
we detail a list of the relevant results of this paper.
\begin{itemize}
\item Theorem~\ref{Theorem:trace-general}
    and Corollary~\ref{corollary:trace}
    are novel Sobolev--trace inequalities.
\item Theorem~\ref{theorem:new-BA} establishes novel 
$\W^{1,\p}$-stable right-inverses of the divergence
satisfying mixed boundary conditions.
\item Theorems~\ref{theorem:main-1}
and~\ref{theorem:main-2} are concerned
with Sobolev--Poincar\'e inequalities
on broken Sobolev spaces.
\item Corollaries~\ref{corollary:main-1} and~\ref{corollary:main-2}
are variants of the Sobolev--Poincar\'e inequalities above,
which are more suitable for certain nonconforming finite element spaces,
including Crouzeix-Raviart spaces;
\item Remark~\ref{remark:applicability} discusses
the applicability of Sobolev-Poincar\'e inequalities
on broken Sobolev spaces to nonconforming Galerkin methods
(including IPDG, WOPSIP, Crouzeix-Raviart)
for linear and nonlinear problems.
\end{itemize}

\noindent We are aware that results
similar to Theorems~\ref{theorem:main-1}
and~\ref{theorem:main-2},
and Corollaries~\ref{corollary:main-1} and~\ref{corollary:main-2}
are work in progress also in~\cite{Cangiani-Dong-Georgoulis-Lin:2025}.
The results discussed therein are, however,
derived undertaking a different avenue,
which is closer to that in~\cite{Brenner:2003}.

\subsection{Domains of interest and Sobolev spaces} 
\label{subsection:domains-Sobolev}
Let~$\Omega$ be a polytopic, open, Lipschitz domain in $\Rbb^d$,
$d$ integer larger than~$1$,
with boundary $\partial\Omega$. Furthermore,
\begin{equation} \label{rho}
\text{$\Omega$ is star-shaped with respect
to a ball $\Brho$ of radius~$\rho$;
the diameter of~$\Omega$ is~$\hOmega$.}
\end{equation}
Assumption~\eqref{rho} can be generalized to domains
as in~\eqref{more-general-domains} below,
i.e., domains that can be decomposed into shape-regular simplicial meshes;
cf. also Remarks~\ref{remark:BA-general-domains}
and~\ref{remark:SP-general-domains},
and Appendix~\ref{appendix:RI-simplices} below.
Even more general domains might be considered,
following e.g. \cite[Sect. 5]{Guzman-Salgado:2021}
or \cite[Remark A.1]{BeiraodaVeiga-Canuto-Nochetto-Vacca:2021}.
We partition the boundary~$\Gamma:= \partial\Omega$
of~$\Omega$ into
\begin{equation} \label{splitting-Gamma}
\Gamma=\Gamma_{\rm D}\cup\Gamma_{\rm N},
\end{equation}
$\GammaD$ having nonzero $(d-1)$ Hausdorff measure
for the sake of simplicity,
and $\Gamma_{\rm D}\cap\Gamma_{\rm N} = \emptyset$.

Given $X \subset \overline{\Omega}$
with diameter~$\h_X$,
we consider Lebesgue spaces of order~$p$
consisting of Lebesgue measurable functions with finite norm
\[
\Norm{v}_{\L^\p(X)}^p
:= \int_X \vert v \vert^p.
\]
Analogously, one defines Lebesgue spaces~$\L^\p(\partial X)$
on~$\partial X$.
$\L^\p_0(X)$ is the space of functions in $\L^\p(X)$
with zero average over~$X$.
We denote the gradient, Laplacian, and divergence operators
by~$\nabla$, $\Delta$, and $\div$;
$\nablau$ and $\Deltabold$ denote
the tensor gradient and the vector Laplacian operators
(these operators act on rows).
More in general, given a multi-index~$\boldalpha$ in~$\Nbb^d$
with size~$s$, $D^{\boldalpha}$ denotes the tensor of dimension~$s$
containing all mixed derivatives of order~$s$.
Given a positive integer~$k$,
and a real number $\p$ in $[1,\infty)$,
$\W^{k,\p}(X)$ denotes the space
of $\L^\p(X)$ functions with
distributional derivatives $D^{\boldalpha}$
of order~$k$ in $\L^\p(X)$.
We introduce the seminorms and norms
\[
\SemiNorm{v}_{\W^{k,\p}(X)}
:= \big(\sum_{\vert \boldalpha \vert=k} 
        \Norm{D^{\boldalpha} v}^{\p}_{\L^\p(X)}\big)^{\frac1\p},
\qquad\qquad
\Norm{v}_{\W^{k,\p}(X)}
:= \Big( \sum_{\ell=0}^k \big(\h_X^{\ell-k}
            \SemiNorm{v}_{\W^{\ell,\p}(X)}\big)^\p \Big)^\frac1\p.
\]
Interpolation theory is used to define
Sobolev spaces of positive noninteger order~$s$.
We shall use the boldface type
to denote vector fields (and the corresponding spaces);
for instance, scalar and vector generic
Lebesgue and Sobolev spaces
are $\L^\p(X)$ and $\Lbf^\p(X)$,
and $\W^{k,p}(X)$ and $\Wbf^{k,p}(X)$.

There exists a bounded linear map, called trace operator
\cite[Section 3.2]{Ern-Guermond:2021},
from $\W^{1,\p}(X)$ to $\L^\p(\partial X)$,
which acts as the restriction to~$\partial X$ for continuous functions.
The image of~$\W^{1,\p}(X)$ is $\W^{1-\frac1\p,\p}(\partial X)$;
see Corollary~\ref{corollary:trace} below
for a clearer statement.
The subspace of~$\W^{1,\p}(X)$ of functions with
zero trace over $\partial X$ is denoted by $\W^{1,\p}_0(X)$.
For all $\widetilde\Gamma \subset \partial X$
with nonzero ($d-1$)-dimensional Hausdorff measure,
we also define
\[
\W^{1,\p}_{\widetilde\Gamma}(X)
:= \{ v \in \W^{1,\p}(X) \mid v_{|\widetilde\Gamma}=0 \} .
\]
A Poincar\'e--Steklov inequality holds true
\cite[Lemma 3.30]{Ern-Guermond:2021}:
for all $\p$ in $[1,\infty)$, there exists
a positive constant~$\CPS{\p}{X}$ such that
\begin{equation} \label{Poincare-Steklov}
\Norm{v}_{\L^\p(X)}
\le \CPS{\p}{X} \h_X \SemiNorm{v}_{\W^{1,\p}(X)}
\qquad\qquad\qquad
\forall v \in \W^{1,\p}_{\Gammatilde}(X),
\quad \Gammatilde \subseteq \partial X.
\end{equation}
A similar inequality holds true
for functions in $\W^{1,\p}(X) \cap L^\p_0(X)$;
with an abuse of notation, we denote the involved constant
with the same symbol.

Continuous embeddings (denoted with the symbol
$\hookrightarrow$) of Sobolev spaces
onto Lebesgue spaces hold true
\cite[Section 2.3.2]{Ern-Guermond:2021}:
given $\ell$ positive and~$\p$ larger than or equal to~$1$,
\begin{itemize}
\item if $\ell \p <d$,
$\W^{\ell,\p}(X) \hookrightarrow \L^{q}(X)$
for all $q$ in $[\p,\frac{\p d}{d-\ell\p}]$;
\item if $\ell \p = d$,
$\W^{\ell,\p}(X) \hookrightarrow \L^{q}(X)$
for all $q$ in $[\p,\infty)$.
\end{itemize}
For future convenience, we spell out the generic
Sobolev embedding estimate:
\begin{equation} \label{Sobolev-embedding}
\Norm{v}_{\L^{q}(X)}
\le \CSob{q}{\ell}{\p}{X} 
    \h_X^{\frac{d}{q}- \frac{d}{\p} + \ell} \Norm{v}_{\W^{\ell,\p}(X)}
\qquad\qquad\qquad
\forall v \in \W^{\ell,\p}(X).
\end{equation}
Given an index~$\p$ in $[1,\infty)$, we define
\begin{equation} \label{indices}
    \pprime :=  \frac{\p}{\p-1};
    \qquad\qquad
    \pstar := 
    \begin{cases}
    \frac{\p d}{d-\p}   & \text{if } \p < d  \\
    \infty              & \text{otherwise};
    \end{cases}
    \qquad\qquad
    \psharp := 
    \begin{cases}
    \frac{\p(d-1)}{d-\p}    & \text{if } \p < d \\
    \infty                  & \text{otherwise}.
    \end{cases}
\end{equation}
The first index in~\eqref{indices}
is the conjugate index of~$\p$;
the second one relates to Sobolev embeddings in dimension~$d$
and $\W^{1,\p}(X)$;
the third one to Sobolev embeddings
on the $(d-1)$-dimensional boundary
and $\W^{1-\frac{1}{\p},\p}(\partial X)$
(combine the trace and Sobolev embedding theorems).

\subsection{Meshes and broken Sobolev spaces} \label{subsection:meshes-broken}
We consider families of meshes~$\{\taun\}$
where each $\taun$ is a finite collection of disjoint, closed,
polytopic elements such that
$\overline{\Omega}=\bigcup_{\E\in\taun}\E$.
Differential operators defined piecewise
over~$\taun$ are denoted with a subscript~$\h$;
for instance, $\nablaboldh$ denotes the broken gradient operator over~$\taun$.
For each $\E$ in $\taun$, $\partial \E$ and $h_\E$ denote
the boundary and the diameter of~$\E$, respectively.

We associate each~$\taun$ with a set $\Fcaln$ covering the mesh skeleton,
i.e., $\bigcup_{\E\in\taun}\partial \E = \bigcup_{F\in\Fcaln} F$.
A facet $\F$ in $\Fcaln$ is a hyperplanar,
closed, and connected subset of $\overline{\Omega}$
with positive $ (d{-}1) $-dimensional
Hausdorff measure such that
\begin{itemize}
\item either there exist distinct $\E_{1,F}$ and $\E_{2,F}$ in $\taun$
such that $F\subseteq\partial \E_{1,F}\cap\partial \E_{2,F}$
and $\F$ is called an \emph{internal facet}
(sometimes also an \emph{interface}),
\item or there exists $\E_F$ in $\taun$ such that
$F\subseteq\partial \E_F\cap\partial\Omega$
and $\F$ is called a \emph{boundary facet}.
\end{itemize}
We collect interfaces and boundary facets
in~$\FcalnI$ and~$\FcalnB$.
For all $\E$ in~$\taun$, $\FcalE$ and $\nbfE$ denote
the set of facets contained in~$\partial\E$ and
the outward unit normal to $\partial\E$.
For a given~$\F$ in~$\Fcaln$,
we fix once and for all one of the two unit normal vectors $\nbfF$.

We assume that the mesh boundary skeleton is compatible
with splitting \eqref{splitting-Gamma},
i.e., $\FcalnB=\FcalnD \cup \FcalnN$ where
\begin{equation} \label{boundary-facets}
\FcalnD:=\{ F\in\FcalnB\,|\, F\subseteq \Gamma_{\rm D}\},
\qquad\qquad\qquad
\FcalnN:=\{ F\in\FcalnB\,|\, F\subseteq \Gamma_{\rm N}\}.
\end{equation}
Following \cite[Assumption 2.1]{Cangiani-Dong-Georgoulis:2016},
we demand the following regularity assumptions:
\begin{itemize}
\item for all~$\E$ in any~$\taun$,
there exists a partition $\mathfrak{T}_\E$
of~$\E$ into non-overlapping $d$-dimensional simplices;
\item there exists a universal,
positive constant~$\gamma$ such that,
for all~$\E$ in any $\taun$ and
all $\T$ in $\frakTE$ with $\partial\T \cap \partial\E \neq \emptyset$,
given~$\FET$ the $(d-1)$-dimensional simplex
$\partial\T \cap \partial\E$,
\begin{equation} \label{regularity-mesh}
\gamma\hE \le d \vert \T \vert \vert \FET \vert^{-1}.
\end{equation}
 \end{itemize}
\begin{remark}[On the regularity of the family of meshes]
\label{remark:regularity-meshes}
The regularity assumptions neither impose a restriction
on the number of facets per element nor on the facets' size.
Moreover, for sequences~$\{\taun\}$
of simplicial meshes,
the above assumptions boil down to the standard shape-regularity assumption.
In particular, classical nonconforming
meshes~\cite{Brenner:2003} are covered by our theory.
\eremk
\end{remark}

Broken Sobolev spaces associated with $\taun$ are defined as 
\[
\W^{1,p}(\taun):=
\left\{u\in L^p(\Omega)\,\mid\, {u}_{|\E}\in \W^{1,p}(\E)
                        \qquad \forall \E\in\taun\right\}.
\]
For every $v\in \W^{1,p}(\taun)$ and $\F$ in $\Fcaln$, the jump operator on $\F$ is given by
$$
\jump{v}_F:=
\begin{cases}
v_{|\E_{1,\F}} \nbf_{\E_{1,\F}} \cdot\nbfF
        + v_{|\E_{2,\F}} \nbf_{\E_{2,\F}} \cdot\nbfF
\quad &\text{if } F \in \FcalnI,\; \F \subset \partial \E_{1,\F}\cap \partial \E_{2,\F},\;\\
v_{|\F} \nbf_{\E_{\F}} \cdot\nbfF
\quad &\text{if } F\in\FcalnB,\; \F \subset \partial\E_\F \cap \partial\Omega.
\end{cases}
$$
With an abuse of notation, we use the same symbols
for vector fields $\vbf$ in~$\Wbf^{1,p}(\taun)$.

\subsection{Main technical tools} \label{subsection:main-technical-tools}
We report here the technical tools that are needed to derive
the results in Section~\ref{subsection:main-results}
below,
and postpone their proofs to
Sections~\ref{section:1st-ingredient}
and~\ref{section:2nd-ingredient} below.

\paragraph*{Sobolev--trace inequalities.}
\begin{theorem}[Sobolev--trace inequalities]
\label{Theorem:trace-general}
Let~$\{\taun\}$ be a family of meshes
as in Section~\ref{subsection:meshes-broken},
$\q$ be in $[1,\infty)$,
$s$ be in $[q,\infty)$,
and~$\gamma$ be as in~\eqref{regularity-mesh}.
For all $\E$ in any~$\taun$
and all $v$ in $\W^{1,\frac{s}{s-\q+1}}(\E) \cap L^s(\E)$,
we have
\begin{equation}\label{eq:improved_cti-general}
\Norm{v}_{\L^{\q}(\partial \E)} 
\le \CTR{\q}{s}{d}{\gamma}
    \Big( \hE^{-\frac1\q\big( 1 - \frac d s (s-\q) \big)} \Norm{v}_{\L^{s}(\E)}
    + \hE^{\frac{1}{\qprime} 
        \big( 1 - \frac d s (s-\q) \big)} 
            \Norm{\nablabold v}_{\Lbf^{\frac{s}{s-\q+1}}(\E)} \Big),
\end{equation}
where,
given $\Gamma(\cdot)$
the Euler's Gamma function,
\begin{equation} \label{eq:trace-constant}
\CTR{\q}{s}{d}{\gamma}
:= \left(
    \frac{d\pi^{\frac{d(s-\q)}{2s}}}{\gamma\Gamma\big( \frac d2 + 1 \big)^{\frac{s-\q}{s}}} 
    + \frac{\q-1}{\gamma} \right) ^\frac1\q.
\end{equation}
\end{theorem}
Two special Sobolev--trace
inequalities are an immediate consequence
of Theorem~\ref{Theorem:trace-general}.
\begin{corollary}[Special Sobolev--trace inequalities]
\label{corollary:trace}
Let~$\{\taun\}$ be a family of meshes
as in Section~\ref{subsection:meshes-broken},
$\q$ be in $[1,\infty)$,
$\CTR{\cdot}{\cdot}{\cdot}{\cdot}$ be as in~\eqref{eq:trace-constant},
the indices~$\qprime$, $\qsharp$, and~$\qstar$
be as in~\eqref{indices},
and~$\gamma$ be as in~\eqref{regularity-mesh}.
For all $\E$ in any~$\taun$
and all $v$ in $\W^{1,\q}(\E)$,
we have
\begin{equation}\label{eq:standard_cti}
\Norm{v}_{\L^{\q}(\partial\E)} 
\le \CTR{\q}{\q}{d}{\gamma}
    \big( \hE^{-\frac1\q} \Norm{v}_{\L^{\q}(\E)}
    + \hE^{\frac1{\qprime}} \Norm{\nablabold v}_{\Lbf^\q(\E)} \big)
\end{equation}
and, if~$q$ further belongs to $[1,d)$,
\begin{equation}\label{eq:improved_cti}
\Norm{v}_{\L^{\qsharp}(\partial\E)} 
\le \CTR{\qsharp}{\qstar}{d}{\gamma} \big(\Norm{v}_{\L^{\qstar}(\E)}
    +\Norm{\nablabold v}_{\Lbf^\q(\E)} \big).
\end{equation}
\end{corollary}

Inequality~\eqref{eq:standard_cti} generalizes
the standard continuous trace inequality,
cf. \cite[Lemma~12.15]{Ern-Guermond:2021} and
\cite[Lemma~1.31]{DiPietro-Droniou:2020},
to the case of families of polytopic meshes
as in~\eqref{regularity-mesh}.

\paragraph*{$\W^{1,\p}$-stable right-inverse
of the divergence satisfying
mixed boundary conditions.}
Babu\v ska--Aziz (BA) inequalities are crucial tools in the analysis of PDEs;
see Appendix~\ref{appendix:homogeneous-BCs-identities-Banach}.
The simplest BA inequality states that there exists
a $\W^{1,\p}$-stable right-inverse of the divergence operator
and typically involves vector fields
with either free boundary conditions
or imposing homogeneous Dirichlet boundary conditions
on the boundary~$\Gamma$
of the domain of interest~$\Omega$.
We shall be using the following result,
whose proof can be found, e.g.,
in \cite{Acosta-Duran:2017, Galdi:2011}.
\begin{lemma}[$\W^{1,\p}$-stable right-inverse of the divergence:
homogeneous boundary conditions]
\label{lemma:Galdi}
Given~$\p$ in~$(1,\infty)$,
there exists a positive constant~$C_{BA}(\p,\Omega)$
such that for all $f$ in $\L^p_0(\Omega)$,
it is possible to construct~$\vbf$ in~$\Wbf^{1,p}_{0}(\Omega)$
satisfying
\begin{equation} \label{BA-fully-homogeneous-BCS}
\div\vbf=f ,
\qquad\qquad\qquad
\SemiNorm{\vbf}_{\Wbf^{1,\p}(\Omega)}
\le C_{BA}(\p,\Omega) \Norm{f}_{\L^\p(\Omega)}.
\end{equation}
As discussed in \cite[Lemma III.3.1]{Galdi:2011},
given~$\hOmega$ and~$\rho$ as in~\eqref{rho},
we have the bound
\begin{equation} \label{Galdi's-bound}
C_{BA}(\p,\Omega)
\le \CG{d}{\p} \Big( \frac{\hOmega}{\rho} \Big)^d  
                \Big( 1+ \frac{\hOmega}{\rho} \Big),
\end{equation}
where $\CG{\cdot}{\cdot}$ is a positive constant
only depending on~$d$ and~$\p$.
As discussed in~\cite{Duran:2012}, an estimate
that is sharper than that in~\eqref{Galdi's-bound}
can be derived in the case $\p$ equal to~2.
\end{lemma}

Recall the splitting~\eqref{splitting-Gamma}
of~$\Gamma$ into
the union of~$\GammaD$ and~$\GammaN$,
assume that both sets have
nonzero $(d-1)$-dimensional measure in~$\Gamma$,
and recall the compatibility
condition~\eqref{boundary-facets}.
We derive BA inequalities
for vector fields with homogeneous boundary conditions
only on~$\GammaN$.
To this aim, we prove a bound
with a constant behaving differently depending
on the convexity of the domain~$\Omega$.
In particular, for nonconvex~$\Omega$,
given~$\GammatildeD$ contained in $\GammaD$,
the regularity of the domain implies the existence
of a domain $\Omegaext$ such that
$(\overline{\Omega}\cap\overline{\Omegaext})^\circ 
= \GammatildeD^\circ$;
see Figure~\ref{figure:extension-nonconvex} below
for an example of~$\Omegaext$.
Given
\begin{equation} \label{Omegaext-properties}
\Omegatilde := (\overline{\Omega} \cup \overline{\Omegaext})^\circ,
\end{equation}
we observe that the regularity properties of~$\Omega$ imply that
\begin{equation} \label{rhotilde}
\text{$\Omegatilde$ is star-shaped with respect
to a ball $\Brhotilde$ of radius~$\rhotilde$;
the diameter of~$\Omegatilde$ is~$\hOmegatilde$.}
\end{equation}
Further define
\begin{equation} \label{rhopartialOmega}
\rhoGamma
:= C_{\Gamma}(\Omega)
    \frac12 \min \Big\{ \h_{\Gamma_j} 
                    = \text{diam($\Gamma_j$)} \mid 
    \Gamma_j \text{ is a (d-1)-dimensional facet of } \Omega \Big\},
\end{equation}
where $C_{\Gamma}(\Omega)$
is constant in $(0,1)$
only depending on $\Omega$
and discussed in the proof of Theorem~\ref{theorem:new-BA}.

\begin{theorem}[$\W^{1,\p}$-stable right-inverse of the divergence:
mixed boundary conditions]
\label{theorem:new-BA}
Given~$\p$ in $(1,\infty)$,
there exists a positive constant~$\CBA{\p}{\GammaN}{\Omega}$
such that, for all $f$ in $\L^\p(\Omega)$,
it is possible to construct~$\vbf$ in~$\Wbf^{1,\p}_{\GammaN}(\Omega)$
satisfying
\begin{equation} \label{BA-partial-homogeneous-BCS}
\div\vbf=f,
\qquad\qquad\qquad
\SemiNorm{\vbf}_{\Wbf^{1,\p}(\Omega)}
\le \CBA{\p}{\GammaN}{\Omega} \Norm{f}_{\L^\p(\Omega)}.
\end{equation}
In particular,
given $\CG{\cdot}{\cdot}$ in~\eqref{Galdi's-bound},
$\hOmega$ and~$\rho$ in~\eqref{rho},
$\hOmegatilde$ and~$\rhotilde$ in~\eqref{rhotilde},
$\Omegaext$ in~\eqref{Omegaext-properties},
and~$\rhoGamma$ in~\eqref{rhopartialOmega},
we have the upper bounds
\begin{equation} \label{constants:Babu-Aziz-mixed}
\CBA{\p}{\GammaN}{\Omega}
\le
\begin{cases}
2^\frac{1}{\p} \CG{d}{\p} \Big( \frac{2\hOmega}{\min(\rho,\rhoGamma)} \Big)^d
    \Big( 1+ \frac{2\hOmega}{\min(\rho,\rhoGamma)} \Big)
    & \text{if } \Omega \text{ is convex} \\
\CG{d}{\p} \Big( \frac{\hOmegatilde}{\rhotilde} \Big)^d  
    \Big( 1+ \frac{\hOmegatilde}{\rhotilde} \Big)
    \left( 1 + \frac{\vert \Omega\vert^{\p-1}}{\vert\Omegaext\vert^{\p-1}} \right)^{\frac1\p}
    & \text{if } \Omega \text{ is not convex}. \\
\end{cases}
\end{equation}
\end{theorem}
In the nonconvex case, the bound on $\CBA{\p}{\GammaN}{\Omega}$
depends heavily on the construction of the extended domain~$\Omegaext$
as detailed above.
In particular, e.g. for domains with re-entrant corners,
the construction of an admissible~$\Omegaext$
as in~\eqref{rhotilde} leads to constants
in~\eqref{constants:Babu-Aziz-mixed} (second line)
that might be much larger than for convex domains.

\begin{remark}[$\W^{1,\p}$-stable right-inverse of the divergence on more general domains]
\label{remark:BA-general-domains}
Lemma~\ref{lemma:Galdi}
and Theorem~\ref{theorem:new-BA}
are proved for star-shaped domains~$\Omega$
satisfying~\eqref{rho}.
Based on Appendix~\ref{appendix:RI-simplices},
the bounds proven in these results are also valid on domains that
can be decomposed into shape-regular simplicial meshes,
cf. Theorem~\ref{theorem:Babu-Aziz-general-domain} below.
Of course, the constants involved in the corresponding inequalities change:
in the homogeneous boundary conditions case,
the constant is now independent of~$\Omega$ and depends on the
shape-regularity parameter,
cf.~\eqref{BA-fully-homogeneous-BCS:general-domains} below;
in the mixed boundary conditions case
the dependence is both on the domain~$\Omega$,
$\GammaN$ as in~\eqref{splitting-Gamma},
and the shape-regularity parameter.
\eremk
\end{remark}

\subsection{Main results} \label{subsection:main-results}
We begin by presenting two Sobolev--Poincar\'e inequalities
for functions in broken $\W^{1,p}$ spaces
and postpone their proofs to Sections~\ref{subsection:main-result-1}
and~\ref{subsection:main-result-2} below.
The first result reads as follows.

\begin{theorem}[First kind Sobolev--Poincar\'e inequalities]
\label{theorem:main-1}
Let~$\CBA{\cdot}{\cdot}{\cdot}$,
$\CTR{\cdot}{\cdot}{\cdot}{\cdot}$,
$\CSob{\cdot}{\cdot}{\cdot}{\cdot}$,
and~$\CPS{\cdot}{\cdot}$
be the constants in~\eqref{BA-partial-homogeneous-BCS},
\eqref{eq:improved_cti-general},
\eqref{Sobolev-embedding},
and~\eqref{Poincare-Steklov},
and~$\hOmega$ be as in~\eqref{rho}.
Consider a family of polytopic meshes $\{\taun\}$
as in Section~\ref{subsection:meshes-broken},
$\p$ in $[1,d)$,
and~$\gamma$ as in~\eqref{regularity-mesh}.
Introduce
\begin{equation} \label{C1-C2}
\begin{split}
C_1 & :=  2^\frac1{(\pstar)'}  \CBA{(\pstar)'}{\GammaN}{\Omega}\,
            \CSob{\pprime}{1}{(\pstar)'}{\Omega}
            (1+ \hOmega^{(\pstar)'}
            \CPS{(\pstar)'}{\Omega}^{(\pstar)'})
            ^\frac1{(\pstar)'},\\
C_2 & :=    \CTR{(\psharp)'}{\pprime}{d}{\gamma}
            \big( \CBA{(\pstar)'}{\GammaN}{\Omega}
                    + C_1  \big).
\end{split}
\end{equation}
Then, we have
\begin{equation} \label{main-result:1}
\Norm{v}_{\L^{\pstar}(\Omega)}
\le C_1        \Norm{\nablaboldh v}_{\Lbf^\p(\Omega)}  \\
 + C_2  \Big(\sum_{\F\in\FcalnI\cup\FcalnD} 
    \Norm{\jump{v}}_{\L^{\psharp}(\F)}^{\p} \Big)^{\frac{1}{\p}} 
\qquad\qquad \forall v \in \W^{1,\p}(\taun,\Omega).
\end{equation}
\end{theorem}
We further have an inequality involving weaker norms.
To this aim, we preliminary introduce
\begin{equation} \label{htildeF}
\htildeF
\quad \text{ either equal to } \hF
\text{ or equal to }\quad
\begin{cases}
    \min(\h_{\E_1} , \h_{\E_2}) & \text{if } \F\in\FcalnI,\; \F \subset \partial\E_1\cap\partial\E_2\\
    \h_\E & \text{if } \F\in\FcalnD,\; \F \in\FcalE.
\end{cases}         
\end{equation}

\begin{theorem}[Second kind Sobolev--Poincar\'e inequalities]
\label{theorem:main-2}
Let~$\CBA{\cdot}{\cdot}{\cdot}$,
$\CTR{\cdot}{\cdot}{\cdot}{\cdot}$,
$\CSob{\cdot}{\cdot}{\cdot}{\cdot}$,
and~$\CPS{\cdot}{\cdot}$
be the constants in~\eqref{BA-partial-homogeneous-BCS},
\eqref{eq:improved_cti-general},
\eqref{Sobolev-embedding},
and~\eqref{Poincare-Steklov},
and~$\hOmega$ be as in~\eqref{rho}.
Consider a family of polytopic meshes $\{\taun\}$
as in Section~\ref{subsection:meshes-broken},
$\p$ in $[1,\infty)$,
and~$\gamma$ as in~\eqref{regularity-mesh}.
Introduce
\begin{equation} \label{C3-C4}
\begin{split}
C_3 & := 2^\frac1\p \CBA{(\p\ostar)'}{\GammaN}{\Omega}\,
        \hOmega^{\frac{d}{\pprime}-\frac{d}{(\p\ostar)'}+1}
         \CSob{\pprime}{1}{(\p\ostar)'}{\Omega} 
        (1+ \hOmega \CPS{(\p\ostar)'}{\Omega}),\\
C_4 & :=  \CBA{(\p\ostar)'}{\GammaN}{\Omega}\,
        \CTR{\pprime}{\pprime}{d}{\gamma}
        \big( 1 + \max_{\E\in\taun} \hE^{\pprime} \big)^{\frac{1}{\pprime}}
        (1+ \hOmega \CPS{\pprime}{\Omega}).
\end{split}
\end{equation}
Then, given~$\htildeF$ as in~\eqref{htildeF}, we have
\begin{equation} \label{main-result:2}
\Norm{v}_{\L^{\postar}(\Omega)}
\le C_3 \Norm{\nablaboldh v}_{\Lbf^\p(\Omega)}
    + C_4 \Big(\sum_{\F\in\FcalnI\cup\FcalnD}
    \htildeF^{1-\p} \Norm{\jump{v}}_{\L^{\p}(\F)}^\p \Big)^\frac1\p
    \qquad\qquad \forall v \in \W^{1,\p}(\taun,\Omega).
\end{equation}
\end{theorem}

\begin{remark}[Norms comparison
and finite dimensional spaces]
\label{remark:norms-finite-dimension}
Theorem~\ref{theorem:main-2}
involves a weaker Sobolev--Poincar\'e inequality
compared to that in Theorem~\ref{theorem:main-1}.
However, if we restrict inequality~\ref{main-result:1}
to functions~$\vh$ in a finite dimensional space
such that Lebesgue inverse inequalities
with explicit constants are available
and the facets of the mesh are uniformly shape-regular,
then~\eqref{main-result:1} can be improved to
\[
\Norm{\vh}_{\L^{\pstar}(\Omega)}
\lesssim \Norm{\nablaboldh \vh}_{\Lbf^\p(\Omega)}  \\
 + \Big(\sum_{\F\in\FcalnI\cup\FcalnD} 
    \hF^{1-\p}
    \Norm{\jump{\vh}}_{\L^{\p}(\F)}^{\p} \Big)^{\frac{1}{\p}} .
\]
\eremk
\end{remark}

We also present two corollaries to
Theorems~\ref{theorem:main-1} and~\ref{theorem:main-2}
that are more useful, e.g.,
for Crouzeix-Raviart type discretizations,
and are in the spirit of~\cite{Mardal-Winther:2006}
and \cite[Remark~1.1]{Brenner:2003}.
Their proofs are given in Sections~\ref{subsection:main-result-3}
and~\ref{subsection:main-result-4} below.

Introduce $\PizFcaln$ mapping functions $\L^1(\Fcaln)$ into
their piecewise average over the facets in~$\Fcaln$:
\[
(\PizFcaln v)_{|\F}
:= \frac{1}{\vert\F\vert} \int_\F v
\qquad\qquad\qquad
\forall v \in \L^1(\F),
\quad \forall \F\in\Fcaln.
\]

\begin{corollary} [First kind Sobolev--Poincar\'e inequalities:
averaged version]
\label{corollary:main-1}
Let~$\CTR{\cdot}{\cdot}{\cdot}{\cdot}$,
$\CSob{\cdot}{\cdot}{\cdot}{\cdot}$,
$\CPS{\cdot}{\cdot}$, and~$C_2$
be the constants in~\eqref{eq:improved_cti-general},
\eqref{Sobolev-embedding},
\eqref{Poincare-Steklov}, and~\eqref{C1-C2},
and~$\hOmega$ be as in~\eqref{rho}.
Consider a family of polytopic meshes $\{\taun\}$
as in Section~\ref{subsection:meshes-broken},
$\p$ in $[1,d)$,
and~$\gamma$ as in~\eqref{regularity-mesh}.
Introduce
\begin{equation} \label{C5-C6}
\begin{split}
C_5
& :=  \max_{\E\in\taun} \big[
            \CSob{\pstar}{1}{\p}{\E}
            (1+\hE \CPS{\p}{\E})  \big] \\
& \qquad      + 2^{1-\frac1\p} C_2
            \Big[\max_{\E\in\taun} (\card (\FcalE) ) \Big]
        \CTR{\psharp}{\pstar}{d}{\gamma}\cdot\\
& \qquad\qquad \cdot\big[\max_{\E\in\taun}
            (1+ \CSob{\pstar}{1}{\p}{\E})
                (1+ \hE \CPS{\pstar}{\E}) \big],
        \qquad C_6 := 2^{1-\frac1\p} C_2 .
\end{split}
\end{equation}
Then, we have
\[
\Norm{v}_{\L^{\pstar}(\Omega)}
 \le C_5 \Norm{\nablaboldh v}_{\Lbf^\p(\Omega)} 
    + C_6  \Big( \sum_{\F\in\FcalnI\cup\FcalnD}
    \Norm{\PizFcaln\jump{v}}_{\L^{\psharp}(\F)}^\p \Big)^\frac1\p
    \qquad\qquad \forall v \in \W^{1,\p}(\taun,\Omega).
\]
\end{corollary}

\begin{corollary} [Second kind Sobolev--Poincar\'e inequalities: averaged version]
\label{corollary:main-2}
Let~$\CTR{\cdot}{\cdot}{\cdot}{\cdot}$,
$\CSob{\cdot}{\cdot}{\cdot}{\cdot}$,
$\CPS{\cdot}{\cdot}$, and~$C_4$
be the constants in~\eqref{eq:improved_cti-general},
\eqref{Sobolev-embedding},
\eqref{Poincare-Steklov}, and~\eqref{C3-C4},
and~$\hOmega$ be as in~\eqref{rho}.
Consider a family of polytopic meshes $\{\taun\}$
as in Section~\ref{subsection:meshes-broken},
$\p$ in $[1,\infty)$,
and~$\gamma$ as in~\eqref{regularity-mesh}.
Given~$\htildeF$ as in~\eqref{htildeF}, introduce
\footnotesize\begin{equation*}
\begin{split}
C_7
& := \max_{\E\in\taun} \big[ 
            \hE^{\frac{d}{\p\ostar}-\frac{d}{\p}+1}
            \CSob{\p\ostar}{1}{\p}{\E}
                        (1+ \hE \CPS{\p}{\E})  \big]\\
& \quad  + 2^{1-\frac1\p} C_4
            \Big[\max_{\E\in\taun} \max_{\F\in\FcalE}
        \Big( \frac{\htildeF}{\hE} \Big)^{-1+\frac1\p}\Big]
        \CTR{\p}{\p}{d}{\gamma}
        \big[ \max_{\E\in\taun} (1+ \hE \CPS{\p}{\E}) \big],
        \qquad C_8  := 2^{1-\frac1\p} C_4 .
\end{split}
\end{equation*}\normalsize
Then, we have
\[
\Norm{v}_{\L^{\p\ostar}(\Omega)}
\le C_7 \Norm{\nablaboldh v}_{\Lbf^\p(\Omega)}
      +  C_8 \Big( \sum_{\F\in\FcalnI\cup\FcalnD}
    \htildeF^{1-\p} \Norm{\PizFcaln\jump{v}}_{\L^\p(\F)}^\p \Big)^\frac1\p 
    \qquad\qquad \forall v \in \W^{1,\p}(\taun,\Omega).
\]
\end{corollary}

\begin{remark}[Sobolev-Poincar\'e inequalities for more general domains.]
\label{remark:SP-general-domains}
Theorems~\ref{theorem:main-1} and~\ref{theorem:main-2},
and Corollaries~\ref{corollary:main-1} and~\ref{corollary:main-2}
are proved for star-shaped domains~$\Omega$
satisfying~\eqref{rho}.
Based on Appendix~\ref{appendix:RI-simplices} below,
the bounds proven in these results are also valid on domains that
can be decomposed into shape-regular simplicial meshes.
Of course, the constants involved in these inequalities
will change and there will be extra dependence
on the shape-regularity parameter while using
the Babu\v ska-Aziz inequality,
cf. Theorem~\ref{theorem:Babu-Aziz-general-domain} below.
\eremk
\end{remark}

\begin{remark}[Sobolev-Poincar\'e inequalities for special families of meshes]
\label{remark:mesh-regularity}
Compared to Theorem~\ref{theorem:main-1},
Corollary~\ref{corollary:main-1}
displays a bound that is not robust with respect to
families of meshes of elements
with number of facets arbitrarily increasing.
This can be seen in the factor
$\max_{\E\in\taun} (\card (\FcalE) )$
appearing in the constant~$C_5$ in~\eqref{C5-C6}.
On the other hand, Corollary~\ref{corollary:main-2}
contains bounds that are robust
with respect to the number of facets of an element
and small facets:
it suffices to take the second option
for~$\htildeF$ in~\eqref{htildeF},
and assume that the mesh is locally quasi-uniform
in the sense that neighbouring elements
have uniformly comparable sizes.
\eremk
\end{remark}

\begin{remark}[Applications to nonconforming methods]
\label{remark:applicability}
Theorem~\ref{theorem:main-2} can be used in the analysis of DG~\cite{Burman-Ern:2008}
and polyDG~\cite{Antonietti-Bonetti-Botti:2025}
methods for nonlinear problems.
Theorem~\ref{theorem:main-1} contains the optimal
Lebesgue exponents; yet, the norm on the right-hand
side of~\eqref{main-result:1}
is stronger than that in~\eqref{main-result:2}
and therefore it would require additional inverse 
estimates in order to be applied in a standard DG and polyDG setting.
This would come at the price of constants depending
on the degree of the method,
which is a major hindrance in the proof of
convergence of the $p$- and $hp$-versions
DG methods under minimal regularity assumptions.
Corollary~\ref{corollary:main-1} can be instead
used for the convergence analysis of Crouzeix-Raviart
schemes with optimal Lebesgue exponents;
cf., e.g., \cite{Girault-Wheeler:2008}.
Furthermore, the fact that the jump terms
appear with an $\L^2$ projection onto constants over facets
allows us to recover the weaker exponent in
standard DG norms as well.
Corollary~\ref{corollary:main-2}
may be used to derive similar results
for Crouzeix-Raviart schemes
and be of inspiration in the analysis
of the WOPSIP method~\cite{Brenner-Sung:2015}
for nonlinear problems.
In view of Remarks~\ref{remark:BA-general-domains}
and~\ref{remark:SP-general-domains},
all the above results may turn out to be particularly useful
on fairly general polytopic domains
and methods based on arbitrary agglomerations
of simplicial elements.
\eremk
\end{remark}

\begin{remark}[The higher-order case]
\label{remark:higher-order}
Theorems~\ref{theorem:main-1} and~\ref{theorem:main-2},
and Corollaries~\ref{corollary:main-1} and~\ref{corollary:main-2}
can be extended to the $\W^{k,p}$ setting
by including further jump terms
involving higher order derivatives in the estimates
using \cite[Appendix~A]{Botti-Mascotto:2025}.
\eremk
\end{remark}

\begin{remark}[Comparison with classical results]
\label{remark:comparison-classical-results}
Theorem~\ref{theorem:main-2}
and Corollary~\ref{corollary:main-2}
generalize the results in \cite{Arnold:1982, Brenner:2003}
in the sense that it suffices to pick $p=2$
and observe that $\ostar$ is larger than~$1$.
Compared to the references above, we admit more general meshes
and exhibit bounds with fully explicit constants
in terms of the properties of the domain~$\Omega$
and the underlying family of meshes.
\eremk
\end{remark}

\paragraph*{Outline of the remainder of the paper.}
In Sections~\ref{section:1st-ingredient} and~\ref{section:2nd-ingredient}
we prove the two technical tools discussed in
Section~\ref{subsection:main-technical-tools}, i.e.,
Sobolev--trace inequalities and
the existence of a $\W^{1,\p}$-stable
right-inverse of the divergence
with mixed boundary conditions.
These two results are instrumental in proving the main results
of the paper, i.e., Theorems~\ref{theorem:main-1}
and~\ref{theorem:main-2},
and Corollaries~\ref{corollary:main-1}
and~\ref{corollary:main-2},
which are the topic of Section~\ref{section:main-results}.

\section{Proof of Sobolev--trace inequalities} \label{section:1st-ingredient}
We prove here Theorem~\ref{Theorem:trace-general}.

\noindent \textbf{Step~$1$.}
We show first an inequality
on a simplex $\T$ in $\frakTE$
as in Section~\ref{subsection:meshes-broken}
for a given $\E$ in any~$\taun$
such that~$\FET:=\partial \T \cap \partial \E\ne \emptyset$.
Let~$\mathbf{P}_{\FET}$ be the vertex of $\T$ opposite to~$\FET$.
We proceed as in \cite[Lemma 12.15]{Ern-Guermond:2021}
and \cite[Lemma 1.31]{DiPietro-Droniou:2020}
by considering the lowest order Raviart--Thomas function
\[
\boldsymbol{\phi}_{\FET}(\mathbf{x})
:=\frac{|\FET|}{d|\T|}(\mathbf{x} - \mathbf{P}_{\FET}).
\]
The normal component of $\boldsymbol{\phi}_{\FET}$
is equal to $1$ on $\FET$ and zero on the other facets of $\T$;
moreover, $\div\boldsymbol{\phi}_{\FET}= |\FET| |\T|^{-1}$.
The divergence theorem implies
\[
\begin{split}
\Norm{v}_{\L^{\q}(\FET)}^{\q} 
&   = \int_{\partial \T} |v|^{\q} \boldsymbol{\phi}_{\FET}\cdot\nbfT
    = \int_\T \div(|v|^{\q} \boldsymbol{\phi}_{\FET})  \\
&   = \int_\T \frac{|\FET|}{|\T|}|v|^{\q} 
        + \int_\T \frac{\q|\FET|}{d|\T|} v |v|^{\q-2} \ 
          \nablabold v \cdot (\mathbf{x} - \mathbf{P}_{\FET})
    =: \mathcal{I}_1 + \mathcal{I}_2.
\end{split}
\]
Using H\"older's inequality
with indices $s/\q$ and $s/(s-\q)$,
and observing that~\eqref{regularity-mesh} implies
$|\FET||\T|^{-1} \le d (\gamma \hE)^{-1}$,
we bound the term $\mathcal{I}_1$ as follows:
\[
\mathcal{I}_1
\le \Norm{v}_{\L^{s}(\T)}^{\q} \frac{|\FET|}{|\T|} |\T|^{\frac{s-\q}{s}} 
\le  d \gamma^{-1} \hE^{-1} |\T|^{\frac{s-\q}{s}} 
            \Norm{v}_{\L^{s}(\T)}^{\q} . 
\]
As for the term $\mathcal{I}_2$,
we remark that $\Vert \mathbf{x} - \mathbf{P}_{\FET} \Vert_{\ell^2} \le \hE$
for all $\mathbf{x}$ in $\T$,
apply H\"older's inequality with
indices~$s/(\q-1)$ and~$s/(s-\q+1)$,
observe that $\qprime(\qsharp -1)=\qstar$,
and use again~\eqref{regularity-mesh} to infer
\[
\begin{split}
\mathcal{I}_2
& \le \frac{\q \hE |\FET|}{d|\T|} 
    \int_\T |v|^{\q -1} \vert \nablabold v \vert
 = \frac{\q \hE |\FET|}{d|\T|}
    \Norm{v}_{\L^{s}(\T)}^{\q-1}
    \Norm{\nablabold v}_{\Lbf^{\frac{s}{s-\q+1}}(\T)}
\le \frac\q\gamma \Norm{v}_{\L^{s}(\T)}^{\q-1}
    \Norm{\nablabold v}_{\Lbf^{\frac{s}{s-\q+1}}(\T)}.
\end{split}
\]
Gathering the estimate of $\mathcal{I}_1$ and $\mathcal{I}_2$,
we obtain a multiplicative trace inequality
for simplicial elements~$\E$ reading  
\begin{equation}\label{eq:cti_simplex-general}
\Norm{v}_{\L^\q(\FET)}^{\q}
\le \gamma^{-1} \Norm{v}_{\L^{s}(\T)}^{\q -1} 
    \left( d \hE^{-1} |\T|^{\frac{s-\q}{s}} \Norm{v}_{\L^{s}(\T)}
        + \q \Norm{\nablabold v}_{\Lbf^{\frac{s}{s-\q+1}}(\T)} \right).
\end{equation}

\noindent\textbf{Step~$2$.}
Let~$\E$ be as in Step~1.
We define
\[
\tildefrakTE
:= \{ \T \in \frakTE \mid
        \FET:= \partial\T\cap\partial\E \ne \emptyset\}.
\]
Using~\eqref{eq:cti_simplex-general},
we get
\[
\begin{aligned}
& \Norm{v}_{\L^{\q}(\partial \E)}^{\q} 
= \sum_{\T \in \tildefrakTE}
    \Norm{v}_{\L^{\q}(\FET)}^{\q}
\le  \frac{d}{\gamma h_{\E}} 
    \sum_{\T \in \tildefrakTE} 
    |\T|^{\frac{s-\q}{s}} 
    \Norm{v}_{\L^{s}(\T)}^{\q}
    + \frac{\q}{\gamma} \sum_{\T \in \tildefrakTE}
        \Norm{v}_{\L^{s}(\T)}^{\q-1}
        \Norm{\nablabold v}_{\Lbf^{\frac{s}{s-\q+1}}(\T)}.
\end{aligned}
\]
H\"older's inequality for sequences
with indices~$s/\q$ and~$s/(s-\q)$,
and~$s/(\q-1)$ and~$s/(s-\q+1)$ implies
\[
\begin{split}
\Norm{v}_{\L^{\q}(\partial \E)}^{\q} 
& \le \frac{d}{\gamma h_\E}
    \Big(\sum_{\T \in \tildefrakTE} 
     |\T|
    \Big)^{\frac{s-\q}{s}}
    \Big(\sum_{\T \in \tildefrakTE} 
    \Norm{v}_{\L^{s}(\T)}^{s} \Big) ^{\frac\q s}\\
& \quad + \frac{\q}{\gamma}
    \Big(\sum_{\T \in \tildefrakTE}
    \Norm{\nablabold v}_{\Lbf^{\frac{s}{s-\q+1}}(\T)}^{\frac{s}{s-\q+1}}
    \Big)^{\frac{s-\q+1}{s}}
    \Big(\sum_{\T \in \tildefrakTE}
    \Norm{v}_{\L^{s}(\T)}^{s} \Big)^{\frac{\q-1}{s}}.
\end{split}
\]
We further have the trivial inclusion
$\bigcup_{\T\in\tildefrakTE} \T \subset \E$.
In light of this, we obtain
\begin{equation} \label{new:discussion}
\Norm{v}_{\L^{\q}(\partial \E)}^{\q}
\le  \frac{d}{\gamma\hE}  |\E|^{\frac{s-\q}{s}} 
        \Norm{v}_{\L^{s}(\E)}^{\q} 
      + \frac{\q}{\gamma} \Norm{\nablabold v}_{\Lbf^{\frac{s}{s-\q+1}}(\E)}
        \Norm{v}_{\L^{s}(\E)}^{\q-1}.
\end{equation}
Let $\Gamma(\cdot)$ denote the Euler's Gamma function.
Given the unit ball~$\mathcal{B}_d$ in $\Rbb^d$,
we have that
\begin{equation} \label{volume:bdd:Euler}
|\E|
\le \hE^d |\mathcal{B}_d|
=\hE^d \frac{\pi^{\frac d2}}{\Gamma \big( \frac d2 + 1 \big)}.
\end{equation}
Furthermore, the Young inequality
$ab\le\frac{a^r}{r}+\frac{b^{r'}}{r'}$
with $r=\q/(\q-1)$ holds true, which implies
\[\Norm{\nablabold v}_{\Lbf^{\frac{s}{s-\q+1}}(\E)}
    \Norm{v}_{\L^{s}(\E)}^{\q-1} 
\le (\q-1)
    \hE^{-1+\frac d s (s-\q)}
    \frac{\Norm{v}_{\L^{s}(\E)}^{\q}}{\q} 
    + \hE^{\frac{\q}{\qprime} \big( 1 - \frac{d}{s}(s-\q)\big)} 
    \frac{\Norm{\nablabold v}_{\Lbf^{\frac{s}{s-\q+1}}(\E)}^{\q}}{\q} .
\]
Combining the three displays above, we get
\small\begin{equation*}
\begin{split}
 \Norm{v}_{\L^{\q}(\partial\E)}^{\q}
 &\le  \left(
    \frac{d\pi^{\frac{d(s-\q)}{2s}}}{\gamma\Gamma\big( \frac d2 + 1 \big)^{\frac{s-\q}{s}}} 
    + \frac{\q-1}{\gamma} \right)
    \hE^{-1+\frac d s (s-\q)}
        \Norm{v}_{\L^{s}(\E)}^\q
    + \hE^{\frac{\q}{\qprime} \big( 1 - \frac{d}{s}(s-\q)\big)} 
        \Norm{\nablabold v}_{\Lbf^{\frac{s}{s-\q+1}}(\E)}^\q \\
& \le  \left(
    \frac{d\pi^{\frac{d(s-\q)}{2s}}}{\gamma\Gamma\big( \frac d2 + 1 \big)^{\frac{s-\q}{s}}} 
    + \frac{\q-1}{\gamma} \right)
    \Big( \hE^{\frac1\q \big(-1+\frac d s (s-\q)\big)}
        \Norm{v}_{\L^{s}(\E)}
    + \hE^{\frac{1}{\qprime} \big( 1 - \frac{d}{s}(s-\q)\big)} 
        \Norm{\nablabold v}_{\Lbf^{\frac{s}{s-\q+1}}(\E)}
        \Big)^\q .
\end{split}
\end{equation*}\normalsize
Inequality~\eqref{eq:improved_cti-general}
follows by taking the power~$\q^{-1}$
on both sides.

\section{$\W^{1,\p}$-stable right-inverse of the divergence with mixed boundary conditions} \label{section:2nd-ingredient}
We prove here Theorem~\ref{theorem:new-BA}.
Let~$\GammatildeN$ be such that $\GammaN \subseteq \GammatildeN$
and
\begin{equation} \label{gammatildeD}
\GammatildeD := \Gamma\setminus\GammatildeN
\end{equation}
is a single $(d-1)$-dimensional facet of~$\Omega$
We have
\begin{equation} \label{bound-BA}
\CBA{\p}{\GammaN}{\Omega}
\le \CBA{\p}{\GammatildeN}{\Omega}.
\end{equation}
In fact,
using that $\W^{1,\p}_{\GammatildeN}(\Omega)$
is contained in $\W^{1,\p}_{\GammaN}(\Omega)$,
a right-inverse of the divergence
in $\W^{1,\p}_{\GammatildeN}(\Omega)$
is also a right-inverse of the divergence
in $\W^{1,\p}_{\GammaN}(\Omega)$.
In light of~\eqref{bound-BA}, it suffices to prove
an upper bound on the constant
appearing on the right-hand side.

\paragraph*{Upper bound on~$\CBA{\p}{\GammatildeN}{\Omega}$ on convex domains.}
Let~$\Omega$ be convex.
Without loss of generality, we assume
that~$\GammatildeD$ lies on the hyperplane $x_1=0$.
In particular, due to the convexity assumption,
it is possible to construct
\[
\Omegaext
:= \{ \xbftilde\in\Rbb^d \mid \xbftilde=(-x_1,x_2,\dots,x_d)
    \; \forall \xbf\in\Omega \}
\]
such that $(\overline\Omega\cap\overline{\Omegaext})^\circ = \GammatildeD^\circ$.
We define $\Omegatilde:= \Omega \cup \Omegaext \cup \GammatildeD^\circ$
and set~$\Gammatilde$ to be its boundary.

The diameter~$\hOmegatilde$ of~$\Omegatilde$ is bounded by $2\hOmega$.
By construction, $\Omegatilde$ is star-shaped with respect
to a ball~$\Brhotilde$ of radius $\rhotilde$,
which is larger than or equal
to the minimum between $\rho$
and~$\rhoGamma$ in~\eqref{rhopartialOmega}
for a suitable choice of the constant~$C_{\Gamma}(\Omega)$
in~\eqref{rhopartialOmega},
which only depends on the convexity
of the domain.
We refer to Figure~\ref{figure:extension-convex}
for a graphical representation of the above construction
in two dimensions (the constant~$C_\Gamma(\Omega)$
is in this example is close to~$1$).

\begin{figure}[ht] 
\begin{center}
\begin{tikzpicture}[scale=1] 
\draw[black, thick] (0,1) -- (-3,2) -- (-3,-2) -- (0,-1);
\draw[green, thick, dashed] (0,1) -- (3,2) -- (3,-2) -- (0,-1);
\draw[green, thick, dashed, opacity=0.4] (0,1) -- (-1.5,.5);
\draw[green, thick, dashed, opacity=0.4] (0,-1) -- (-1.5,-.5);
\draw[black, thick, dashed, opacity=0.4] (0,1) -- (1.5,.5);
\draw[black, thick, dashed, opacity=0.4] (0,-1) -- (1.5,-.5);
\draw (-3,1.7) node[black, below, right] {$\Omega$};
\draw (2.95,-1.45) node[green, below, left] {$\Omegaext$};
\draw[orange, thick] (0,-1) -- (0,1);
\draw (0,0.2) node[black, left] {$\orange{\GammatildeD}$};
\draw[blue, thick] (0,0) circle (.9cm); 
\draw (-1,1) node[black] {$\blue{\Brhotilde}$};
\draw (3,1.5) node[black, left] {$\red{\hOmegatilde}$};
\draw[red, thick] (-3,-2) -- (3,2);
\draw (0,1.5) node[black] {$\Omegatilde$};
\end{tikzpicture}
\end{center}
\caption{Extended domain~$\Omegatilde$
(with black and dashed green boundary) for a convex domain~$\Omega$.
The opaque dashed lines are used to determine
the radius of the ball~$\Brhotilde$.}
\label{figure:extension-convex}
\end{figure}

Given~$f$ in $\L^\p(\Omega)$, we define
$\ftilde$ in $\L^\p_0(\Omegatilde)$ as
\begin{equation} \label{ftilde:convex}
\ftilde(\xbf)
:= \begin{cases}
f(x_1,x_2,\dots,x_d)   & \text{if } \xbf=(x_1,x_2,\dots,x_d) \in \Omega    \\
-f(-x_1,x_2,\dots,x_d)  & \text{if } \xbf=(x_1,x_2,\dots,x_d) \in \Omegaext .
\end{cases}
\end{equation}
Let~$\vbftilde$ in $\Wbf^{1,\p}_0(\Omegatilde)$
be the right-inverse of the divergence
applied to~$\ftilde$ in $\L^\p_0(\Omegatilde)$.
Defining $\vbf = \vbftilde_{|\Omega}$
in $\Wbf^{1,\p}_{\GammatildeN}(\Omega)$
and using~\eqref{BA-fully-homogeneous-BCS}, we deduce
\begin{equation} \label{birichino}
\SemiNorm{\vbf}_{\Wbf^{1,\p}(\Omega)}
= \SemiNorm{\vbftilde}_{\Wbf^{1,\p}(\Omega)}
\le \SemiNorm{\vbftilde}_{\Wbf^{1,\p}(\Omegatilde)}
\le \CBA{\p}{\Gammatilde}{\Omegatilde} \Norm{\ftilde}_{\L^\p(\Omegatilde)}
\le 2^\frac1\p \CBA{\p}{\Gammatilde}{\Omegatilde} \Norm{f}_{\L^\p(\Omegatilde)}.
\end{equation}
This entails
\[
\CBA{\p}{\GammaN}{\Omega}
\overset{\eqref{bound-BA}}{\le}
    \CBA{\p}{\GammatildeN}{\Omega}
\overset{\eqref{birichino}}{\le}
    2^\frac1\p \CBA{\p}{\Gammatilde}{\Omegatilde}.
\]
Further using Lemma~\ref{lemma:Galdi}
and recalling that~$\rhoGamma$
is defined in~\eqref{rhopartialOmega}, we write
\[
\CBA{\p}{\Gammatilde}{\Omegatilde}
\le \CG{d}{\p} \Big( \frac{\hOmegatilde}{\rhotilde} \Big)^d  
                \Big( 1+ \frac{\hOmegatilde}{\rhotilde} \Big)
\le \CG{d}{\p} \Big( \frac{2\hOmega} 
                {\min(\rho,\rhoGamma)} \Big)^d  
                \Big( 1+ \frac{2\hOmega}{\min(\rho,\rhoGamma)} \Big).
\]

\paragraph*{Upper bound on~$\CBA{\GammatildeN}{\p}{\Omega}$ on nonconvex domains.}
Let $\Omegaext$ as detailed in Section~\ref{subsection:main-technical-tools};
see Figure~\ref{figure:extension-nonconvex}
for a graphical representation of that construction.
\begin{figure}[ht] 
\begin{center}
\begin{tikzpicture}[scale=1] 
\draw[black, thick] (-3,2) -- (-2,-1) -- (2,-1) -- (3,2) -- (0,0.5);
\draw[green, thick, dashed] (-3,2) -- (0,2) -- (0,0.5);
\draw (2.8,1.55) node[black, below, left] {$\Omega$};
\draw (0,1.75) node[green, below, left] {$\Omegaext$};
\draw[orange, thick] (-3,2) -- (0,0.5);
\draw (-1.9,1.7) node[orange] {$\GammatildeD$};
\draw[magenta, thick, dashed] (0,0.5) -- (0,-1);
\draw[magenta, thick, dashed] (0,0.5) -- (-15/7,-4/7);
\draw[blue, thick] (-.6,-0.45) circle (0.51cm); 
\draw (-.6,-.45) node[blue] {$\Brhotilde$};
\draw (1.9,-.9) node[black, above] {$\red{\hOmega}$};
\draw[red, thick] (-3,2) -- (2,-1);
\draw (0.3,1) node[black] {$\Omegatilde$};
\end{tikzpicture}
\end{center}
\caption{Extended domain~$\Omegatilde$
(with black and dashed green boundary)
for a nonconvex~$\Omega$.}
\label{figure:extension-nonconvex}
\end{figure}
Given $f$ in~$\L^\p(\Omega)$
and~$\fbarOmega$ its average over~$\Omega$,
define $\ftilde$ in~$\L^\p_0(\Omegatilde)$ as
\begin{equation} \label{ftilde:nonconvex}
\ftilde(\xbf)
:= \begin{cases}
    f(\xbf) & \text{if } \xbf \in \Omega \\
    -\fbarOmega \frac{\vert \Omega \vert}{\vert \Omegaext \vert}
            & \text{if } \xbf \in \Omegaext.
\end{cases}
\end{equation}
We have the following bounds:
\[
\begin{split}
\Norm{\ftilde}_{\L^\p(\Omegatilde)}^\p
& \le  \Norm{f}_{\L^\p(\Omega)}^\p
      + \Norm{\fbarOmega \vert\Omega\vert \vert\Omegaext\vert^{-1}}_{\L^\p(\Omegaext)}^\p 
  \le \Norm{f}_{\L^\p(\Omega)}^\p
    + \frac{\vert \Omega\vert^{\p-1}}{\vert\Omegaext\vert^{\p-1}}
        \Norm{f}_{\L^\p(\Omega)}^\p .
\end{split}
\]
Let~$\vbftilde$ in $\Wbf^{1,\p}_0(\Omegatilde)$
be the right-inverse of the divergence
applied to~$\ftilde$ in $\L^\p_0(\Omegatilde)$.
Defining $\vbf = \vbftilde_{|\Omega}$
in $\Wbf^{1,\p}_{\GammatildeN}(\Omega)$
and using~\eqref{BA-fully-homogeneous-BCS}, we deduce
\begin{equation} \label{brighella}
\begin{split}
\SemiNorm{\vbf}_{\Wbf^{1,\p}(\Omega)}
= \SemiNorm{\vbftilde}_{\Wbf^{1,\p}(\Omega)}
\le \SemiNorm{\vbftilde}_{\Wbf^{1,\p}(\Omegatilde)}
& \le \CBA{\p}{\Gammatilde}{\Omegatilde} \Norm{\ftilde}_{\L^\p(\Omegatilde)}\\
& \le \CBA{\p}{\Gammatilde}{\Omegatilde} 
    \left( 1 + \frac{\vert \Omega\vert^{\p-1}}{\vert\Omegaext\vert^{\p-1}} \right)^{\frac1\p}
    \Norm{f}_{\L^\p(\Omega)}.
\end{split}
\end{equation}
We deduce the upper bound
\[
\CBA{\p}{\GammaN}{\Omega}
\overset{\eqref{bound-BA}}{\le}
    \CBA{\p}{\GammatildeN}{\Omega}
\overset{\eqref{brighella}}{\le}
    \CBA{\p}{\Gammatilde}{\Omegatilde} 
    \left( 1 + \frac{\vert \Omega\vert^{\p-1}}{\vert\Omegaext\vert^{\p-1}} 
    \right)^{\frac1\p}.
\]
Further using Lemma~\ref{lemma:Galdi},
given~$\hOmegatilde$ and~$\rhotilde$
as in~\eqref{rhotilde},
we write
\[
\CBA{\p}{\Gammatilde}{\Omegatilde} 
\le \CG{d}{\p} \Big( \frac{\hOmegatilde}{\rhotilde} \Big)^d  
                \Big( 1+ \frac{\hOmegatilde}{\rhotilde} \Big).
\]

\begin{remark}[On the smoothness of the right-inverse of the
divergence on nonconvex domains.]
\label{remark:smoothness-RI}
The construction for nonconvex domains directly applies
to the convex case for any extended domain~$\Omegaext$;
in particular, this applies also to the case
that $\GammatildeD$
in~\eqref{gammatildeD} is only
part of a facet of the domain~$\Omega$.
The price to pay is that~$\ftilde$ as in~\eqref{ftilde:convex}
is continuous across~$\GammatildeN$,
whereas~$\ftilde$ as in~\eqref{ftilde:nonconvex} is not.
Therefore, the approach for convex domains
yields smoother right-inverses of the divergence
with available stability estimates
as discussed in \cite[Appendix~A]{Botti-Mascotto:2025}.
\eremk
\end{remark}

\section{Proof of the main results} \label{section:main-results}
We prove here the main results of the paper,
namely Theorems~\ref{theorem:main-1} and~\ref{theorem:main-2},
and Corollaries~\ref{corollary:main-1} and~\ref{corollary:main-2}
in Sections~\ref{subsection:main-result-1},
\ref{subsection:main-result-2},
\ref{subsection:main-result-3},
and~\ref{subsection:main-result-4}.

\subsection{Proof of the first main result} \label{subsection:main-result-1}
We prove here Theorem~\ref{theorem:main-1}.
Given~$\p$ in $[1,d)$, let $s=(\pstar)'$.
By duality, we can write
\[
\Norm{v}_{\L^{\pstar}(\Omega)}
= \sup_{\varphi \in L^s(\Omega)} 
        \frac{ _{\L^{\pstar}(\Omega)}(v,\varphi)_{\L^s(\Omega)}}{\Norm{\varphi}_{\L^s(\Omega)}}.
\]
From Theorem~\ref{theorem:new-BA},
we know that there exists~$\CBA{s}{\GammaN}{\Omega}$ such that
for all~$\varphi$ in~$\L^s(\Omega)$ it is possible to construct~$\sigmabold$
in~$\Wbf^{1,s}_{\GammaN}(\Omega)$ satisfying
\[
\div\sigmabold = \varphi ,
\qquad\qquad\qquad
\SemiNorm{\sigmabold}_{\Wbf^{1,s}(\Omega)}
\le \CBA{s}{\GammaN}{\Omega} \Norm{\varphi}_{\L^s(\Omega)}.
\]
We combine the two displays above and readily deduce
\small\begin{equation} \label{1st:IBP}
\begin{split}
 \Norm{v}_{\L^{\pstar}(\Omega)}
 &\le \CBA{s}{\GammaN}{\Omega}\sup_{\sigmabold \in \Wbf^{1,s}_{\GammaN}(\Omega)} 
    \frac{_{\L^{\pstar}(\Omega)}(v,\div\sigmabold)_{\L^s(\Omega)}}{\SemiNorm{\sigmabold}_{\Wbf^{1,s}(\Omega)}} \\
& = \CBA{s}{\GammaN}{\Omega} \sup_{\sigmabold \in \Wbf^{1,s}_{\GammaN}(\Omega)} 
    \frac{_{\Lbf^{\p}(\Omega)}(-\nablaboldh v,\sigmabold)_{\Lbf^{\pprime}(\Omega)}
      + \sum_{\F\in\Fcaln}{}
      _{\L^{\psharp}(\F)}(\jump{v},\sigmabold\cdot\nbfF)_{\L^{(\psharp)'}(\F)}}
      {\SemiNorm{\sigmabold}_{\Wbf^{1,s}(\Omega)}}.
\end{split}
\end{equation}\normalsize
We check that the two terms in the numerator
on the right-hand side above are well-defined.
\begin{itemize}
\item Since~$v$ belongs to~$\W^{1,p}_{\GammaN}(\taun,\Omega)$,
we clearly have that~$\nablaboldh v$ is in~$\Lbf^\p(\Omega)$.
\item Recall that $\sigmabold$ belongs to $\Wbf^{1,s}(\Omega)$
with $s=(\pstar)'$;
the Sobolev embedding
$\Wbf^{1,s}(\Omega) \hookrightarrow \Lbf^{q}(\Omega)$
holds true~\eqref{Sobolev-embedding}
for all $q$ in $[1,\sstar]$, where we recall
$\sstar=(ds)/(d-s)$.
To conclude, we have to show that~$\sstar=\pprime$.
Standard manipulations imply
\[
\begin{split}
\sstar
& = \frac{s d}{d-s}
  = \frac{\pstar d}{d(\pstar-1)-\pstar}
  = \frac{\pstar d}{\pstar(d-1)-d}
  = \frac{\p d}{\p d-\p-d+\p}
  = \frac{\p}{\p-1}
  = \pprime.
\end{split}
\]
\item The trace of~$v$ in~$\W^{1,p}_{\GammaN}(\taun,\Omega)$
belongs to $\W^{1-\frac1\p,\p}(\F)$ for all facets~$\F$ in~$\Fcaln$;
see, e.g., \cite[Chapter~3]{Ern-Guermond:2021}.
Using the Sobolev embedding (in dimension $d-1$)
$\W^{1-\frac1\p,\p}(\F) \hookrightarrow L^{q}(\F)$
for all $q$ such that
\[
q \le \frac{\p(d-1)}{(d-1)-(\p-1)}
= \frac{\p(d-1)}{d-\p}
= \psharp.
\]
In particular, the trace of $v$ belongs to~$\L^{\psharp}(\F)$
for all facets~$\F$.
\item The trace of~$\sigmabold$ in~$\Wbf^{1,\s}(\Omega)$
belongs to $\Wbf^{1-\frac1\s,\s}(\F)$
for all facets~$\F$ in~$\Fcaln$.
Using the Sobolev embedding~\eqref{Sobolev-embedding}
(in dimension $d-1$)
$\Wbf^{1-\frac1\s,\s}(\F) \hookrightarrow \Lbf^{q}(\F)$
for all $r$ such that,
proceeding as above,
$r \le \ssharp.$
We are left with proving $\ssharp=(\psharp)'$.
Standard manipulations imply
\[
\begin{split}
\ssharp
& = \frac{\s(d-1)}{d-\s}
  = \frac{\left(\frac{\p d}{d-\p}\right)'(d-1)}{d-\left(\frac{\p d}{d-\p}\right)'}
= \frac{\frac{\p d}{\p d - d +\p}(d-1)}{d-\frac{\p d}{\p d - d+\p}}
  = \frac{\p(d-1)}{\p d -d}
  = \frac{\psharp}{\psharp-1}
  = (\psharp)' .
\end{split}
\]
\end{itemize}
We now estimate the two terms in the numerator
on the right-hand side of~\eqref{1st:IBP}.
The Sobolev embedding~\eqref{Sobolev-embedding}
$\Wbf^{1,s}(\Omega) \hookrightarrow \Lbf^{\pprime}(\Omega)$
and the Poincar\'e--Steklov inequality~\eqref{Poincare-Steklov}
imply
\begin{equation} \label{1st:1st-term}
\begin{split}
\Norm{\sigmabold}_{\Lbf^{\pprime}(\Omega)}
& \le \hOmega^{\frac{d}{\pprime}-\frac{d}{s}+1}
    \CSob{\pprime}{1}{s}{\Omega} \Norm{\sigmabold}_{\Wbf^{1,s}(\Omega)} \\
& \le \hOmega^{\frac{d}{\pprime}-\frac{d}{s}+1}
    \CSob{\pprime}{1}{s}{\Omega} 
    (1+\hOmega^s \CPS{s}{\Omega})^s)^\frac1s
         \SemiNorm{\sigmabold}_{\Wbf^{1,s}(\Omega)}.
\end{split}
\end{equation}
Moreover,
using two H\"older's inequalities
with indices $\psharp$ and $\ssharp$,
and $\p$ and $\pprime=\sstar= \ssharp \ostar$
(for sequences),
and Jensen's inequality for sequences,
we get
\small\[
\begin{split}
& \sum_{\F\in\Fcaln}{}
      _{\L^{\psharp}(\F)}(\jump{v},\sigmabold\cdot\nbfF)_{\L^{(\psharp)'}(\F)}
= \sum_{\F\in\FcalnI \cup \FcalnD}{}
      _{\L^{\psharp}(\F)}(\jump{v},\sigmabold\cdot\nbfF)_{\L^{(\psharp)'}(\F)}\\
& \le \sum_{\F\in\FcalnI \cup \FcalnD}{}
    \Norm{\jump{v}}_{\L^{\psharp}(\F)}
    \Norm{\sigmabold\cdot\nbfF}_{\L^{\ssharp}(\F)} 
  \le \Big( \sum_{\F\in\FcalnI \cup \FcalnD}
        \Norm{\jump{v}}_{\L^{\psharp}(\F)}^\p
        \Big)^{\frac1\p}
    \Big( \sum_{\E\in\taun} \sum_{\F\in\FcalE}
        \Norm{\sigmabold\cdot\nbfF}_{\L^{\ssharp}(\F)}^{\sstar}
        \Big)^{\frac1{\sstar}} \\
& \le \Big( \sum_{\F\in\FcalnI \cup \FcalnD}
        \Norm{\jump{v}}_{\L^{\psharp}(\F)}^\p
        \Big)^{\frac1\p}
    \Big( \sum_{\E\in\taun} \big(\sum_{\F\in\FcalE}
        \Norm{\sigmabold\cdot\nbfF}_{\L^{\ssharp}(\F)}^{\ssharp} \big)^{\ostar}
        \Big)^{\frac{1}{\sstar}} \\
& = \Big( \sum_{\F\in\FcalnI \cup \FcalnD}
        \Norm{\jump{v}}_{\L^{\psharp}(\F)}^\p
        \Big)^{\frac1\p}
    \Big( \sum_{\E\in\taun}
        \Norm{\sigmabold\cdot\nbfE}_{\L^{\ssharp}(\partial\E)}^{\sstar}
        \Big)^{\frac{1}{\sstar}} .
\end{split}
\]\normalsize   
We apply the Sobolev--trace inequality~\eqref{eq:improved_cti}
to the last term on the right-hand side
and Jensen's inequality for sequences,
and get
\[
\begin{split}
& \sum_{\F\in\Fcaln}{}
      _{\L^{\psharp}(\F)}(\jump{v},\sigmabold\cdot\nbfF)_{\L^{(\psharp)'}(\F)} \\
& \le 2^\frac1s \CTR{\ssharp}{\sstar}{d}{\gamma}
    \Big( \sum_{\F\in\FcalnI \cup \FcalnD}
        \Norm{\jump{v}}_{\L^{\psharp}(\F)}^\p
        \Big)^{\frac1\p}
        \Big( \sum_{\E\in\taun}
        (\Norm{\sigmabold}_{\Lbf^{\sstar}(\E)}^{\sstar}
        + \SemiNorm{\sigmabold}_{\Wbf^{1,s}(\E)}^{\sstar}) 
        \Big)^{\frac{1}{\sstar}} \\
& \le 2^\frac1s \CTR{\ssharp}{\sstar}{d}{\gamma}
    \Big( \sum_{\F\in\FcalnI \cup \FcalnD}
        \Norm{\jump{v}}_{\L^{\psharp}(\F)}^\p
        \Big)^{\frac1\p}
        \Big(
        \Norm{\sigmabold}_{\Lbf^{\sstar}(\Omega)}^{\sstar}
        + (\sum_{\E\in\taun} \SemiNorm{\sigmabold}_{\Wbf^{1,s}(\E)}^s)^{\frac{\sstar}{s}}
        \Big)^{\frac{1}{\sstar}} \\
& \le  2^\frac1s
    \CTR{\ssharp}{\sstar}{d}{\gamma}
    \Big( \sum_{\F\in\FcalnI \cup \FcalnD}
        \Norm{\jump{v}}_{\L^{\psharp}(\F)}^\p
        \Big)^{\frac1\p}
        \Big(
        \Norm{\sigmabold}_{\Lbf^{\sstar}(\Omega)}
        + \SemiNorm{\sigmabold}_{\Wbf^{1,s}(\Omega)}
        \Big).
\end{split}
\]\normalsize
We apply a Sobolev embedding as in~\eqref{Sobolev-embedding}
and the Poincar\'e--Steklov
inequality in~\eqref{Poincare-Steklov}
to get
\begin{equation} \label{1st:2nd-term}
\begin{split}
& \sum_{\F\in\Fcaln}{}
      _{\L^{\psharp}(\F)}(\jump{v},\sigmabold\cdot\nbfF)_{\L^{(\psharp)'}(\F)} \\
& \le 2^\frac1s \CTR{\ssharp}{\sstar}{d}{\gamma}
    \big( 1+    \hOmega^{\frac{d}{\sstar}-\frac{d}{s}+1}
                \CSob{\sstar}{1}{s}{\Omega}
            (1+ \hOmega^s \CPS{s}{\Omega})^s)^\frac1s
        \big)\cdot\\
& \qquad\qquad\qquad\qquad
    \cdot \Big( \sum_{\F\in\FcalnI \cup \FcalnD}
        \Norm{\jump{v}}_{\L^{\psharp}(\F)}^\p
        \Big)^{\frac1\p}
        \SemiNorm{\sigmabold}_{\Wbf^{1,s}(\Omega)}.
\end{split}
\end{equation}
The assertion follows
recalling $\s = (\pstar)'$,
combining~\eqref{1st:IBP}, \eqref{1st:1st-term}, and~\eqref{1st:2nd-term},
and noting that $\frac{d}{\pprime}-\frac{d}{(\pstar)'}+1$ equals $0$.

\subsection{Proof of the second main result} \label{subsection:main-result-2}
We prove here Theorem~\ref{theorem:main-2}.
Throughout, we employ the notation
in the proof of Theorem~\ref{theorem:main-1}.
Introduce
\begin{equation} \label{definition:s-bis}
s
:= (\postar)'
= \left(\frac{d \p}{d-1}\right)'
= \frac{d\p}{d\p-d+1}.
\end{equation}
By duality, we write
\[
\Norm{v}_{\L^{\postar}(\Omega)}
= \sup_{\varphi \in L^s(\Omega)} 
        \frac{ _{\L^{\postar}(\Omega)}(v,\varphi)_{\L^s(\Omega)}}{\Norm{\varphi}_{\L^s(\Omega)}}.
\]
From Theorem~\ref{theorem:new-BA},
we know that there exists~$\CBA{s}{\GammaN}{\Omega}$ such that
for all~$\varphi$ in~$\L^s(\Omega)$
it is possible to construct~$\sigmabold$ in~$\Wbf^{1,s}_{\GammaN}(\Omega)$ satisfying
\[
\div\sigmabold = \varphi ,
\qquad\qquad\qquad
\SemiNorm{\sigmabold}_{\Wbf^{1,s}(\Omega)}
\le \CBA{s}{\GammaN}{\Omega} \Norm{\varphi}_{\L^s(\Omega)}.
\]
We combine the two displays above and get
\begin{equation} \label{1st:IBP-bis}
\begin{split}
& \Norm{v}_{\L^{\p \ostar}(\Omega)}
 \le \CBA{s}{\GammaN}{\Omega} \sup_{\sigmabold \in \Wbf_{\GammaN}^{1,s}(\Omega)} 
    \frac{_{\L^{\postar}(\Omega)}(v,\div\sigmabold)_{\L^s(\Omega)}}{\SemiNorm{\sigmabold}_{\Wbf^{1,s}(\Omega)}} \\
& = \CBA{s}{\GammaN}{\Omega} \sup_{\sigmabold \in \Wbf_{\GammaN}^{1,s}(\Omega)} 
    \frac{_{\Lbf^{\p}(\Omega)}(-\nablaboldh v,\sigmabold)_{\Lbf^{\pprime}(\Omega)}
      + \sum_{\F\in\Fcaln}{}
      _{\L^{\p}(\F)}(\jump{v},\sigmabold\cdot\nbfF)_{\L^{\pprime}(\F)}}
      {\SemiNorm{\sigmabold}_{\Wbf^{1,s}(\Omega)}}.
\end{split}
\end{equation}
We check that the inner products in the numerator
on the right-hand side above are well-defined.
\begin{itemize}
\item Since~$v$ belongs to~$\W^{1,p}(\taun,\Omega)$,
we clearly have that~$\nablaboldh v$ is in~$\Lbf^\p(\Omega)$.
\item Recall that $\sigmabold$
belongs to $\Wbf^{1,s}_{\GammaN}(\Omega)$
with $s$ as in~\eqref{definition:s-bis};
the Sobolev embedding
$\Wbf^{1,s}(\Omega) \hookrightarrow \Lbf^{q}(\Omega)$
holds true~\eqref{Sobolev-embedding}
for all $q$ in $[1,\sstar]$, where we recall
$\sstar=(ds)/(d-s)$.
To conclude, we show that~$\sstar=\pprime \ostar$,
which is larger than~$\pprime$.
Standard manipulations imply
\[
\begin{split}
\sstar
& \overset{\eqref{definition:s-bis}}{=}
    \frac{(\p\ostar)' d}{d-(\p\ostar)'}
    = \frac{d^2 \p}{d^2\p-d^2+d-d\p} 
    = \frac{d \p}{(d-1)(\p-1)}
    = d'\pprime
    =\ostar\pprime.
\end{split}
\]
\item The trace of functions in~$\W^{1,\p}(\Omega)$
belongs to $\L^\p(\F)$ for all facets~$\F$ in~$\Fcaln$.
\item The trace of~$\Wbf^{1,\s}(\Omega)$
belongs to $\Wbf^{1-\frac1\s,\s}(\F)$
for all facets~$\F$ in~$\Fcaln$.
Using the Sobolev embedding~\eqref{Sobolev-embedding}
(in dimension $d-1$)
$\Wbf^{1-\frac1\s,\s}(\F) \hookrightarrow \Lbf^{q}(\F)$
for all $q$ such that
\[
\begin{split}
q
& \le \frac{s(d-1)}{(d-1)-(s-1)}
    \overset{\eqref{definition:s-bis}}{=}
    \frac{(d-1)d\p}{d^2\p-d^2+d-d\p}
    = \frac{(d-1)\p}{(d-1)(\p-1)}
    = \frac{\p}{\p-1} = \pprime.
\end{split}
\]
\end{itemize}
As for the first term on the right-hand side of~\eqref{1st:IBP-bis},
we exploit the Sobolev embedding~\eqref{Sobolev-embedding}
and the Poincar\'e--Steklov inequality~\eqref{Poincare-Steklov},
and arrive at
\begin{equation} \label{1st:1st-term-bis}
\begin{split}
& _{\Lbf^{\p}(\Omega)}(-\nablaboldh v,\sigmabold)_{\Lbf^{\pprime}(\Omega)}
   \le \Norm{\nablaboldh v}_{\Lbf^\p (\Omega)}
        \Norm{\sigmabold}_{\Lbf^{\pprime}(\Omega)} \\
&   \le \hOmega^{\frac{d}{\pprime}-\frac{d}{(\p\ostar)'}+1}
        \CSob{\pprime}{1}{(\p\ostar)'}{\Omega} 
        \big(1+ \hOmega \CPS{(\p\ostar)'}{\Omega})\big)
        \Norm{\nablaboldh v}_{\Lbf^\p (\Omega)}
        \SemiNorm{\sigmabold}_{\Wbf^{1,s}(\Omega)} .
\end{split}
\end{equation}
Recall that~$\htildeF$ is defined in~\eqref{htildeF}.
As for the second term on the right-hand side of~\eqref{1st:IBP-bis},
using two H\"older's inequalities
with indices $\p$ and $\pprime$,
and $\p$ and $\pprime$ (for sequences), we get
\begin{equation} \label{2nd:2nd-term-a}
\begin{split}
& \sum_{\F \in \FcalnI\cup\FcalnD}{}
    _{\L^{\p}(\F)}(\jump{v},\sigmabold\cdot\nbfF)_{\L^{\pprime}(\F)}
    \le \sum_{\F \in \FcalnI\cup\FcalnD}{}
        \htildeF^{\frac{1-\p}{\p}} \Norm{\jump{v}}_{\L^\p(\F)}
        \htildeF^{-\frac{1-\p}{\p}} \Norm{\sigmabold\cdot\nbfF}_{\L^{\pprime}(\F)}\\
& \le \Big( \sum_{\F \in \FcalnI\cup\FcalnD}{} 
            \htildeF^{1-\p} \Norm{\jump{v}}_{\L^\p(\F)}^\p \Big)^\frac{1}{\p}
    \Big( \sum_{\F \in \FcalnI\cup\FcalnD}{} 
            \htildeF \Norm{\sigmabold\cdot\nbfF}_{\L^{\pprime}(\F)}^{\pprime}    
                        \Big)^\frac{1}{\pprime} .
\end{split}
\end{equation}
We cope with the second term on the right-hand side of~\eqref{2nd:2nd-term-a}
separately.
Using that $\htildeF$ is smaller than or equal to~$\hE$
for $\E$ such that $\F$ belongs to~$\FcalE$,
the trace inequality~\eqref{eq:standard_cti},
Jensen's inequality for sequences,
the Poincar\'e--Steklov inequality~\eqref{Poincare-Steklov},
and standard manipulations yields
\small\[
\begin{split}
& \Big( \sum_{\F \in \FcalnI\cup\FcalnD}{} 
            \htildeF \Norm{\sigmabold\cdot\nbfF}_{\L^{\pprime}(\F)}^{\pprime}    
                        \Big)^\frac{1}{\pprime}
    \le \Big( \sum_{\E \in \taun}{}  \hE
            \sum_{\F\in\FcalE}
            \Norm{\sigmabold\cdot\nbfF}_{\L^{\pprime}(\F)}^{\pprime}    
                        \Big)^\frac{1}{\pprime}\\
& = \Big( \sum_{\E \in \taun}{}  \hE
            \Norm{\sigmabold\cdot\nbfF}_{\L^{\pprime}(\partial\E)}^{\pprime}
            \Big)^\frac{1}{\pprime}
    \le 2^\frac1\p
    \CTR{\pprime}{\pprime}{d}{\gamma}
    \Big( \sum_{\E \in \taun}{}  \hE
            \Big( \hE^{-1} \Norm{\sigmabold}_{\L^{\pprime}(\E)}^{\pprime}
                + \hE^{\frac{\pprime}{\p}}
                    \SemiNorm{\sigmabold}_{\Wbf^{1,\pprime}(\E)}^{\pprime}
                    \Big)  \Big)^\frac{1}{\pprime}\\
& = 2^\frac1\p \CTR{\pprime}{\pprime}{d}{\gamma}
    \Big( \sum_{\E \in \taun}{}
             \Big( \Norm{\sigmabold}_{\L^{\pprime}(\E)}^{\pprime}
                + \hE^{\pprime}
                    \SemiNorm{\sigmabold}_{\Wbf^{1,\pprime}(\E)}^{\pprime}
            \Big)\Big)^\frac{1}{\pprime}\\
& \le 2^\frac1\p \CTR{\pprime}{\pprime}{d}{\gamma}
    \big( 1 + \max_{\E\in\taun} \hE^{\pprime} \big)^{\frac{1}{\pprime}}
    \Big( \Norm{\sigmabold}_{\L^{\pprime}(\Omega)}^{\pprime}
        +\SemiNorm{\sigmabold}_{\Wbf^{1,\pprime}(\Omega)}^{\pprime}
        \Big)^\frac{1}{\pprime}\\
&   \le 2^\frac1\p \CTR{\pprime}{\pprime}{d}{\gamma}
    \big( 1 + \max_{\E\in\taun} \hE^{\pprime} \big)^{\frac{1}{\pprime}}
    (1+ \hOmega \CPS{\pprime}{\Omega})
    \SemiNorm{\sigmabold}_{\Wbf^{1,\pprime}(\Omega)}.
\end{split}
\]\normalsize
A combination of this display with~\eqref{2nd:2nd-term-a} gives
\small\begin{equation} \label{2nd:2nd-term}
\begin{split}
& \sum_{\F \in \FcalnI\cup\FcalnD}{}
    _{\L^{\p}(\F)}(\jump{v},\sigmabold\cdot\nbfF)_{\L^{\pprime}(\F)} \\
& \le 2^\frac1\p \CTR{\pprime}{\pprime}{d}{\gamma}
    \big( 1 + \max_{\E\in\taun} \hE^{\pprime} \big)^{\frac{1}{\pprime}}
    (1+ \hOmega \CPS{\pprime}{\Omega})
    \SemiNorm{\sigmabold}_{\Wbf^{1,\pprime}(\Omega)}
    \Big( \sum_{\F \in \FcalnI\cup\FcalnD}{}
            \htildeF^{1-\p} \Norm{\jump{v}}_{\L^\p(\F)}^\p \Big)^\frac{1}{\p}.
\end{split}    
\end{equation}\normalsize
The assertion follows combining \eqref{1st:IBP-bis},
\eqref{1st:1st-term-bis},
and~\eqref{2nd:2nd-term}.

\subsection{Proof of the third main result} \label{subsection:main-result-3}
We prove here Corollary~\ref{corollary:main-1}.
Here and in the following section,
let $\Piztaun$ be the piecewise average operator over~$\taun$ as
\[
(\Piztaun v)_{|\E}
:= \frac{1}{\vert\E\vert} \int_\E v
\qquad\qquad\qquad
\forall v \in \L^1(\E),
\quad \forall \E\in\taun.
\]
The triangle inequality gives
\begin{equation} \label{eq:cor1-0}
\Norm{v}_{\L^{\pstar}(\Omega)}
\le \Norm{v-\Piztaun v}_{\L^{\pstar}(\Omega)} 
    + \Norm{\Piztaun v}_{\L^{\pstar}(\Omega)}.
\end{equation}
We focus on the first term on the right-hand side:
using the Sobolev embedding~\eqref{Sobolev-embedding} and
the Poincar\'e-Steklov inequality~\eqref{Poincare-Steklov} yields
\[
\Norm{v-\Piztaun v}_{\L^{\pstar}(\Omega)} 
    \le \max_{\E\in\taun} \big( 
                \hE^{\frac{d}{\pstar}-\frac{d}{\p}+1}
                \CSob{\pstar}{1}{\p}{\E}
                (1+ \hE \CPS{\p}{\E})  \big)
            \Norm{\nablaboldh v}_{\Lbf^\p(\Omega)}.
\]
Observing that  $\frac{d}{\pstar}-\frac{d}{\p}+1$ equals $0$
and combining the two displays above entail
\begin{equation} \label{eq:cor1-1}
\Norm{v}_{\L^{\pstar}(\Omega)}
\le \max_{\E\in\taun} \big( 
            \CSob{\pstar}{1}{\p}{\E}
            (1+ \hE \CPS{\p}{\E})  \big) \Norm{\nablaboldh v}_{\Lbf^\p(\Omega)}
    + \Norm{\Piztaun v}_{\L^{\pstar}(\Omega)}.
\end{equation}
We are left with estimating the second term on the right-hand side.
We apply Theorem~\ref{theorem:main-1} and get,
for~$C_2$ as in~\eqref{C1-C2},
\[
\Norm{\Piztaun v}_{\L^{\pstar}(\Omega)}
\le C_2 \Big( \sum_{\F\in\FcalnI\cup\FcalnD}
        \Norm{\jump{\Piztaun v}}_{\L^{\psharp}(\F)}^\p \Big)^\frac1\p .
\]
Recall that $\PizFcaln$ is the piecewise average operator over~$\Fcaln$.
We estimate each jump term separately:
\[
\Norm{\jump{\Piztaun v}}_{\L^{\psharp}(\F)}
\le  \Norm{\jump{\Piztaun v} - \PizFcaln \jump{v}}_{\L^{\psharp}(\F)}
        +  \Norm{\PizFcaln \jump{v}}_{\L^{\psharp}(\F)}.
\]
A direct computation reveals that
\begin{equation}
\label{eq:annabalci}
\Norm{\jump{\Piztaun v} - \PizFcaln \jump{v}}_{\L^{\psharp}(\F)}
= \Norm{\PizFcaln \jump{v-\Piztaun v}}_{\L^{\psharp}(\F)}
\le \Norm{\jump{v-\Piztaun v}}_{\L^{\psharp}(\F)}.
\end{equation}
We combine the three displays above and get
\begin{equation} \label{eq:cor1-2}
\begin{split}
\Norm{\Piztaun v}_{\L^{\pstar}(\Omega)}
& \le 2^{1-\frac1\p} C_2 \Big( \sum_{\F\in\FcalnI\cup\FcalnD}
    \Norm{\PizFcaln\jump{v}}_{\L^{\psharp}(\F)}^\p \Big)^\frac1\p \\
&  \qquad + 2^{1-\frac1\p} C_2 \Big( \sum_{\F\in\FcalnI\cup\FcalnD}
        \Norm{\jump{v-\Piztaun v}}_{\L^{\psharp}(\F)}^\p \Big)^\frac1\p .
\end{split}
\end{equation}
We are left to handle the second term on the right-hand side,
as the first one is fine as it is.
H\"older's inequality for sequences
with indices $\psharp/\p$ and $\psharp/(\psharp-\p)$
reveals that
\[
\begin{split}
& \Big( \sum_{\F\in\FcalnI\cup\FcalnD}
    \Norm{\jump{v-\Piztaun v}}_{\L^{\psharp}(\F)}^\p \Big)^\frac1\p
    \le \Big( \sum_{\E\in\taun}  \sum_{\F \in \FcalE}
        \Norm{v-\Piztaun v}_{\L^{\psharp}(\F)}^\p \Big)^\frac1\p \\
& \le \Big[\max_{\E\in\taun} (\card (\FcalE) )
        ^{\frac{\psharp-\p}{\p}} \Big]
        \Big( \sum_{\E\in\taun} 
            \Norm{v-\Piztaun v}_{\L^{\psharp}(\partial\E)}^{\p}
            \Big)^\frac1\p .
\end{split}
\]
Further using the trace inequality~\eqref{eq:improved_cti},
the Poincar\'e-Steklov inequality~\eqref{Poincare-Steklov},
and the Sobolev embedding~\eqref{Sobolev-embedding},
and recalling that $\frac{d}{\pstar}-\frac{d}{\p}+1$ equals~$0$,
we deduce
\footnotesize\begin{equation} \label{eq:cor1-3}
\begin{split}
& \Big( \sum_{\F\in\FcalnI\cup\FcalnD}
    \Norm{\jump{v-\Piztaun v}}_{\L^{\psharp}(\F)}^\p \Big)^\frac1\p\\
&   \le \Big[\max_{\E\in\taun} (\card (\FcalE) ) \Big]
        \CTR{\psharp}{\pstar}{d}{\gamma} 
        \big[\max_{\E\in\taun} (1+
                \CSob{\pstar}{1}{\p}{\E})
                (1+ \hE \CPS{\pstar}{\E}) \big]
        \Norm{\nablaboldh v}_{\Lbf^\p(\Omega)}.
\end{split}
\end{equation}\normalsize
We combine~\eqref{eq:cor1-0}, \eqref{eq:cor1-1},
\eqref{eq:cor1-2}, and~\eqref{eq:cor1-3}, and write
\[
\begin{split}
\Norm{v}_{\L^{\pstar}(\Omega)}
& \le \Big\{  \max_{\E\in\taun} \big[ 
            \CSob{\pstar}{1}{\p}{\E}
            (1+\hE \CPS{\p}{\E})  \big] \\
& \qquad      + 2^{1-\frac1\p} C_2
            \Big[\max_{\E\in\taun} (\card (\FcalE) ) 
        \CTR{\psharp}{\pstar}{d}{\gamma}\cdot\\
& \qquad\qquad \cdot\big[\max_{\E\in\taun} 
            (1+ \CSob{\pstar}{1}{\p}{\E})
                (1+ \hE \CPS{\pstar}{\E}) \big]
        \Big\} \Norm{\nablaboldh v}_{\Lbf^\p(\Omega)} \\
& \quad + 2^{1-\frac1\p} C_2
        \Big( \sum_{\F\in\FcalnI\cup\FcalnD}
        \Norm{\PizFcaln\jump{v}}_{\L^{\psharp}(\F)}^\p \Big)^\frac1\p ,
\end{split}
\]
which is the assertion.

\subsection{Proof of the fourth main result} \label{subsection:main-result-4}
We prove here Corollary~\ref{corollary:main-2}.
Recall that $\Piztaun$ denotes the piecewise average
operator over~$\taun$.
The triangle inequality gives
\begin{equation} \label{eq:cor2-0}
\Norm{v}_{\L^{\p\ostar}(\Omega)}
\le \Norm{v-\Piztaun v}_{\L^{\p\ostar}(\Omega)} 
    + \Norm{\Piztaun v}_{\L^{\p\ostar}(\Omega)}.
\end{equation}
We focus on the first term on the right-hand side:
using the Sobolev embedding~\eqref{Sobolev-embedding} and
the Poincar\'e-Steklov inequality~\eqref{Poincare-Steklov} yields
\[
\Norm{v-\Piztaun v}_{\L^{\p\ostar}(\Omega)} 
    \le \max_{\E\in\taun} \big(
                \hE^{\frac{d}{\p\ostar}-\frac{d}{\p}+1}
                \CSob{\p\ostar}{1}{\p}{\E}
                        (1+\hE\CPS{\p}{\E})  \big)
            \Norm{\nablaboldh v}_{\Lbf^\p(\Omega)}.
\]
Combining the two displays above entails
\small\begin{equation} \label{eq:cor2-1}
\Norm{v}_{\L^{\p\ostar}(\Omega)}
\le \max_{\E\in\taun} \big( 
            \hE^{\frac{d}{\p\ostar}-\frac{d}{\p}+1}
            \CSob{\p\ostar}{1}{\p}{\E}
            (1+\hE\CPS{\p}{\E})  \big) \Norm{\nablaboldh v}_{\Lbf^\p(\Omega)}
    + \Norm{\Piztaun v}_{\L^{\p\ostar}(\Omega)}.
\end{equation} \normalsize
We are left with estimating the second term on the right-hand side.
For~$C_4$ as in~\eqref{C3-C4}
and~$\htildeF$ as in~\eqref{htildeF},
we apply Theorem~\ref{theorem:main-2} and get
\[
\Norm{\Piztaun v}_{\L^{\p\ostar}(\Omega)}
\le C_4 \Big( \sum_{\F\in\FcalnI\cup\FcalnD}
        \htildeF^{1-\p} \Norm{\jump{\Piztaun v}}_{\L^\p(\F)}^\p \Big)^\frac1\p .
\]
Recall that $\PizFcaln$ is the piecewise average operator over~$\Fcaln$.
Proceeding as in \eqref{eq:annabalci} reveals
\[
\begin{aligned}
& \Norm{\jump{\Piztaun v}}_{\L^\p(\F)}
\le  \Norm{\jump{\Piztaun v} - \PizFcaln \jump{v}}_{\L^\p(\F)}
        +  \Norm{\PizFcaln \jump{v}}_{\L^\p(\F)}\\
&= \Norm{\PizFcaln \jump{v-\Piztaun v}}_{\L^\p(\F)}  +  \Norm{\PizFcaln \jump{v}}_{\L^\p(\F)}
\le \Norm{\jump{v-\Piztaun v}}_{\L^\p(\F)}
+  \Norm{\PizFcaln \jump{v}}_{\L^\p(\F)}.
\end{aligned}
\]
We combine the two displays above and get
\begin{equation} \label{eq:cor2-2}
\begin{split}
\Norm{\Piztaun v}_{\L^{\p\ostar}(\Omega)}
& \le 2^{1-\frac1\p} C_4 \Big( \sum_{\F\in\FcalnI\cup\FcalnD}
    \htildeF^{1-\p} \Norm{\PizFcaln\jump{v}}_{\L^\p(\F)}^\p \Big)^\frac1\p \\
& \quad  + 2^{1-\frac1\p} C_4 \Big( \sum_{\F\in\FcalnI\cup\FcalnD}
        \htildeF^{1-\p} \Norm{\jump{v-\Piztaun v}}_{\L^\p(\F)}^\p \Big)^\frac1\p .
\end{split}
\end{equation}
We are left to handle the second term on the right-hand side,
as the first one is fine as it is:
\[
\begin{split}
& \Big( \sum_{\F\in\FcalnI\cup\FcalnD}
    \htildeF^{1-\p} \Norm{\jump{v-\Piztaun v}}_{\L^\p(\F)}^\p \Big)^\frac1\p
    \le \Big( \sum_{\E\in\taun}  \sum_{\F \in \Fcaln} \htildeF^{1-\p} 
            \Norm{v-\Piztaun v}_{\L^\p(\F)}^\p \Big)^\frac1\p \\
& \le \Big[\max_{\E\in\taun} \max_{\F\in\FcalE}
        \Big( \frac{\htildeF}{\hE} \Big)^{-1+\frac1\p}\Big]
        \Big( \sum_{\E\in\taun}  \hE^{1-\p} 
            \Norm{v-\Piztaun v}_{\L^\p(\partial\E)}^\p \Big)^\frac1\p .
\end{split}
\]
Further using the trace inequality~\eqref{eq:standard_cti}
and the Poincar\'e-Steklov inequality~\eqref{Poincare-Steklov},
we deduce
\small\begin{equation} \label{eq:cor2-3}
\begin{split}
& \Big( \sum_{\F\in\FcalnI\cup\FcalnD}
    \htildeF^{1-\p} \Norm{\jump{v-\Piztaun v}}_{\L^\p(\F)}^\p \Big)^\frac1\p \\
&   \le \Big[\max_{\E\in\taun} \max_{\F\in\FcalE}
        \Big( \frac{\htildeF}{\hE} \Big)^{-1+\frac1\p} \Big]
        \CTR{\p}{\p}{d}{\gamma} 
        \big[ \max_{\E\in\taun} (1+ \hE\CPS{\p}{\E}) \big]
        \Norm{\nablaboldh v}_{\Lbf^\p(\Omega)}.
\end{split}
\end{equation}\normalsize
We combine~\eqref{eq:cor2-0}, \eqref{eq:cor2-1},
\eqref{eq:cor2-2}, and~\eqref{eq:cor2-3}, and write
\footnotesize\[
\begin{split}
\Norm{v}_{\L^{\p\ostar}(\Omega)}
& \le \Big\{ \max_{\E\in\taun} \big[ 
                    \hE^{\frac{d}{\p\ostar}-\frac{d}{\p}+1}
                    \CSob{\p\ostar}{1}{\p}{\E}
                        (1+ \hE \CPS{\p}{\E})  \big]\\
& \qquad      + 2^{1-\frac1\p} C_4
            \Big[\max_{\E\in\taun} \max_{\F\in\FcalE}
        \Big( \frac{\htildeF}{\hE} \Big)^{-1+\frac1\p}\Big]
        \CTR{\p}{\p}{d}{\gamma} 
        \big[ \max_{\E\in\taun} (1+ \hE \CPS{\p}{\E}) \big]
        \Big\} \Norm{\nablaboldh v}_{\Lbf^\p(\Omega)} \\
& \quad + 2^{1-\frac1\p} C_4 
        \Big( \sum_{\F\in\FcalnI\cup\FcalnD}
    \htildeF^{1-\p} \Norm{\PizFcaln\jump{v}}_{\L^\p(\F)}^\p \Big)^\frac1\p ,
\end{split}
\]\normalsize
which is the assertion.

\paragraph*{Acknowledgments.}
The Authors are grateful to J\'er\^ome Droniou
and Martin Vohral\'ik
for useful feedbacks on the first version of the manuscript.
MB and LM have been partially funded by the
European Union (ERC, NEMESIS, project number 101115663);
views and opinions expressed are however those
of the authors only and do not necessarily reflect
those of the EU or the ERC Executive Agency.
LM has been partially funded by MUR (PRIN2022 research grant n. 202292JW3F).
The authors are also members of the Gruppo Nazionale Calcolo Scientifico-Istituto Nazionale di Alta Matematica (GNCS-INdAM).

{\footnotesize
\bibliography{bibliography.bib}}
\bibliographystyle{plain}

\appendix

\section{Inequalities involving the right-inverse of the divergence}
\label{appendix:homogeneous-BCs-identities-Banach}

We recall several inequalities:
given $p$ in $(1,\infty)$,
\begin{itemize}
    \item ({\bf lowest order inf-sup condition}) there exists a positive constant~$\betaz{\p}{\Gamma}{\Omega}$ such that
    \begin{equation} \label{inf-sup-homogeneous-BCS-Banach}
        \inf_{q\in \Lpprimez}
        \sup_{\vbf \in \Wopzd}
        \frac{_{\L^{\p}}(\div \vbf, q)_{\L^{\pprime}}}{\SemiNorm{\vbf}_{\Wbf^{1,\p}(\Omega)} \Norm{q}_{\Lpprime}}
        \ge \betaz{\p}{\Gamma}{\Omega};
    \end{equation}
    \item ({\bf lowest order Ne\v cas-Lions inequality}) there exists a positive constant~$\CNL{\p}{\Gamma}{\Omega}$ such that
    \begin{equation} \label{NL-homogeneous-BCS-Banach}
        \Norm{q}_{\Lpprime} 
        \le \CNL{\p}{\Gamma}{\Omega} \Norm{\nablabold q}_{(\Wbf^{1,\p}_0(\Omega))'}
        \qquad\qquad\qquad \forall q\in \Lpprimez;
    \end{equation}
    \item ({\bf lowest order Babu\v ska--Aziz inequality}) there exists a positive constant~$C_{BA}(\p,\Omega)$ such that
    for all~$q$ in $\L^{p}_0(\Omega)$ it is possible to
    construct $\vbf$ in $\Wopzd$ satisfying
    \begin{equation} \label{BA-homogeneous-BCS-Banach}
        \div\vbf=q,
        \qquad\qquad\qquad
        \SemiNorm{\vbf}_{\Wbf^{1,p}(\Omega)}
        \le C_{BA}(\p,\Omega) \Norm{q}_{\Lp}.
    \end{equation}
\end{itemize}
In the remainder of this section,
cf. the related works~\cite{Horgan-Payne:1983, Costabel-Dauge:2015, Bernardi-Costabel-Dauge-Girault:2016},
we show
\begin{equation} \label{homogeneous-BCs-identities-Banach}
    \betaz{\p}{\Gamma}{\Omega}^{-1} 
    = \CNL{\p}{\Gamma}{\Omega}
    = \CBA{\p}{\Gamma}{\Omega}.
\end{equation}

\paragraph*{Showing $\betaz{\p}{\Gamma}{\Omega}^{-1}=\CNL{\p}{\Gamma}{\Omega}$.}
For all~$q$ in $\Lpprimez$,
an integration by parts and the definition of negative norms
yields
\[
\sup_{\vbf \in \Wopzd}
    \frac{_{\L^\p}(\div \vbf, q)_{\L^{\pprime}}}
    {\SemiNorm{\vbf}_{\Wbf^{1,\p}(\Omega)} \Norm{q}_{\Lpprime}}
= \sup_{\vbf \in \Wopzd}
    \frac{_{(\Wbf_0^{1,\p})'} \langle \nablabold q, \vbf \rangle_{\Wbf^{1,\p}}}
    {\SemiNorm{\vbf}_{\Wbf^{1,\p}(\Omega)} \Norm{q}_{\Lpprime}}
=  \frac{\Norm{\nablabold q}_{(\Wbf^{1,\p}_0(\Omega))'}}{\Norm{q}_{\Lpprime}}.
\]
Taking the infimum over all possible~$q$ entails
\begin{equation} \label{beta=CNLmo}
\betaz{\p}{\Gamma}{\Omega}^{-1}=\CNL{\p}{\Gamma}{\Omega}.
\end{equation}

\paragraph*{Nonlinear operators.}

Let~$\Wbf^{-1,\pprime}(\Omega)$ denote
the dual space of $\Wbf^{1,\p}_0(\Omega)$. Given
\[
\begin{split}
\Vbfpz 
& := \{ \vbf \in \Wbf^{1,\p}_0(\Omega) \mid \div\vbf=0 \}, 
\qquad
\Vbfpperp:=
\{ \vbf \in \Wbf^{1,\p}_0(\Omega)
    \mid (\nablau\vbf,\nablau\wbf)=0
    \quad \forall \wbf \in \Vbfz \},
\end{split}
\]
the following orthogonal decomposition holds true:
\begin{equation} \label{orthogonal-Banach}
\Wbf^{1,\p}_0(\Omega) = \Vbfpperp \oplus^{\Wbf^{1,\p}} \Vbfpz.
\end{equation}
We introduce the operators
\begin{align}
&\Gcalp :\Wbf^{1,\p}_0(\Omega) \to \Lbb^{\pprime}(\Omega),
& \Gcalp(\ubf) = \vert\nablau \ubf\vert^{\p-2} \nablau \ubf, \label{Gcalp}\\
&\Lcalp :\Wbf^{1,\p}_0(\Omega) \to \Wbf^{-1,\pprime}(\Omega) ,
& \Lcalp(\ubf) = \divbf\Gcalp(\ubf), \label{Lcalp}\\
&\Dcalp :\Wbf^{1,\p}_0(\Omega) \to \L^{\pprime}(\Omega),
& \Dcalp(\ubf) = \vert\div \ubf\vert^{\p-2} \div \ubf \label{Dcalp}.
\end{align}
The Banach version of \cite[Lemma 3.43]{John:2016}
gives
\begin{equation} \label{iso-div-Banach}
\text{the $\Dcalp$ operator is a bijection
between $\Vbfpperp$ and $\L^{\pprime}_0(\Omega)$.}
\end{equation}
The Banach version of \cite[Lemmas 3.41 and 3.43]{John:2016}
gives
\begin{equation} \label{iso-nabla-Banach}
\text{the $\nablabold$ operator is a bijection
between $\L^{\pprime}_0(\Omega)$ and $(\Vbfpperp)'$.}
\end{equation}
Noting that~$\Lcalp$ is the solution map
of the $\p$-Laplacian problem with homogeneous
Dirichlet boundary conditions, we also have
\begin{equation} \label{iso-Delta-Banach}
\text{the $\Lcalp$ operator is a bijection
between $\Vbfpperp$ and $(\Vbfpperp)'$.}
\end{equation}

\paragraph*{Nonlinear Cosserat's operator.}
Consider the following nonlinear Stokes' eigenvalue problem:
find $\ubf$ in~$\Wbf^{1,\p}_0(\Omega)$, $\pi$ in $\L^{\pprime}_0(\Omega)$,
and~$\lambda$ in~$\Rbb$ such that
\begin{equation} \label{Stokes-eigenvalue-Banach}
\begin{cases}
-\Lcalp(\ubf) + \nablabold \pi = 0  
    & \text{in } \Omega \\
\div \ubf = \lambda^{\pprime-1} \vert\pi\vert^{\pprime-2}\pi 
    & \text{in } \Omega.
\end{cases}
\end{equation}
In a standard weak formulation,
the solution $(\ubf,\pi)$ belongs to
$\Vbfpperp\times\L^{\pprime}_0(\Omega)$.
Multiplying the first equation by~$\lambda$,
and applying $f(\xi) = \vert\xi\vert^{\p-2}\xi$
and taking the gradient of the second equation give
\[
\begin{cases}
-\lambda\Lcalp(\ubf) + \lambda\nablabold \pi = 0 
    & \text{in } \Omega \\
\nablabold \Dcalp (\ubf) = \lambda \nablabold \pi         
    & \text{in } \Omega.
\end{cases}
\]
We obtain
\[
\nablabold \Dcalp(\ubf)
= \lambda \Lcalp(\ubf)
\qquad \text{ in } (\Vbfpperp)'.
\]
Using~\eqref{iso-Delta-Banach}, we write
\[
\Lcalp^{-1} (\nablabold \Dcalp(\ubf)) 
    = \lambda^{\frac{1}{\p-1}} \ubf
    \quad\text{ in }\quad \Vbfpperp.
\]
Applying~$\Dcalp$ on both sides of the previous identity
and replacing $\Dcalp (\ubf)$ by~$q$ in $\L^{\pprime}_0(\Omega)$,
we get the equivalent nonlinear eigenvalue problem
\[
\Dcalp (\Lcalp^{-1}(\nablabold q)) 
    = \lambda q \quad\text{ in }\quad
    \L^{\pprime}_0(\Omega).
\]
By virtue of~\eqref{iso-div-Banach}, \eqref{iso-nabla-Banach},
and~\eqref{iso-Delta-Banach},
we define the nonlinear Cosserat's operators
(named after the Cosserat brothers \cite{Cosserat-Cosserat:1896})
\begin{equation} \label{Cosserat-operators-Banach}
\Scalp:= \Lcalp^{-1} \nablabold \Dcalp:
    \Vbfpperp \to \Vbfpperp,
\qquad\qquad
\Scalpstar:= \Dcalp \Lcalp^{-1} \nablabold: \L^{\pprime}_0(\Omega) \to \L^{\pprime}_0(\Omega).
\end{equation}
The eigenvalues~$\lambda$ in~\eqref{Stokes-eigenvalue-Banach}
are also eigenvalues of the following nonlinear problems:
\begin{align}
\Lcalp\Scalp \ubf
= \lambda\Lcalp \ubf 
    & \quad\text{ in }\quad (\Vbfpperp)', \label{cosserat.eigen1} \\
|\Scalpstar q|^{\pprime-2}\Scalpstar q
= \lambda^{\pprime-1} |q|^{\pprime-2}q
    & \quad\text{ in }\quad \L^{\p}_0(\Omega) . \label{cosserat.eigen2}
\end{align}
As such, given
\begin{equation} \label{minimum-Cosserat-Banach}
\sigma(\Scalp,\p,\Gamma,\Omega)
:= \text{min. eigenvalue in~\eqref{cosserat.eigen1}},
\quad
\sigma(\Scalpstar,\p,\Gamma,\Omega)
:= \text{min. eigenvalue in~\eqref{cosserat.eigen2}},
\end{equation}
we have
\begin{equation} \label{sigmaSp=sigmaSpstar}
\sigma(\Scalp,\p,\Gamma,\Omega)^{\pprime-1}
= \sigma(\Scalpstar,\p,\Gamma,\Omega).
\end{equation}

\paragraph*{Showing $\sigma(\Scalp,\p,\Gamma,\Omega) = \CBA{\p}{\Gamma}{\Omega}^{-\p}$.}
We characterize~$\sigma(\Scalp,\p,\Gamma,\Omega)$ in~\eqref{minimum-Cosserat-Banach}
by means of the nonlinear Rayleigh quotient:
\begin{equation} \label{sigmaSp=CBA-p}
\begin{split}
\sigma(\Scalp,\p,\Gamma,\Omega)
&   \overset{\eqref{cosserat.eigen1}}{=}
    \inf_{\vbf\in\Vbfpperp}  
    \frac{_{\Lbb^{\pprime}}\langle\Gcalp\Scalp\vbf, \nablau \vbf \rangle_{\Lbb^\p}}{\SemiNorm{\vbf}_{\Wbf^{1,\p}(\Omega)}^\p}
    =
    \inf_{\vbf\in\Vbfpperp} \frac{_{\Wbf^{-1,\pprime}}\langle -\divbf \Gcalp \Lcalp^{-1} \nablabold \Dcalp \vbf, \vbf \rangle_{\Wbf^{1,\p}_0(\Omega)}}{\SemiNorm{\vbf}_{\W^{1,\p}}^\p} \\
&   \overset{\eqref{Lcalp}}{=}
    \inf_{\vbf\in\Vbfpperp} \frac{_{\Wbf^{-1,\pprime}}\langle -\Lcalp\Lcalp^{-1} \nablabold \Dcalp \vbf, \vbf \rangle_{\Wbf^{1,\p}_0}}{\SemiNorm{\vbf}_{\Wbf^{1,\p}(\Omega)}^\p}
    =
    \inf_{\vbf\in\Vbfpperp} \frac{_{\Wbf^{-1,\pprime}}\langle -\nablabold \Dcalp \vbf, \vbf \rangle_{\Wbf^{1,\p}_0}}{\SemiNorm{\vbf}_{\Wbf^{1,\p}(\Omega)}^\p} \\
&   =
    \inf_{\vbf\in\Vbfpperp} \frac{_{\L^{\pprime}}\langle\Dcalp \vbf, \div \vbf\rangle_{\L^\p  }}{\SemiNorm{\vbf}_{\Wbf^{1,\p}(\Omega)}^\p}
    =
    \inf_{\vbf\in\Vbfpperp} \frac{\Norm{\div \vbf}_{\Lbf^\p(\Omega)}^\p}{\SemiNorm{\vbf}_{\Wbf^{1,\p}(\Omega)}^\p} 
    \overset{\eqref{BA-homogeneous-BCS-Banach}}{=:}
    \CBA{\p}{\Gamma}{\Omega}^{-\p}.
\end{split}
\end{equation}

\paragraph*{A remark.}
Let~$q$ be in~$\L^{\pprime}_0(\Omega)$
and $\zetabold$ be in~$\Vbfpperp$,
cf. \eqref{iso-div-Banach} and \eqref{iso-Delta-Banach},
be the solution to
\begin{equation} \label{elliptic-auxiliary-Banach}
\begin{cases}
- \Lcalp \zetabold = \nablabold q  & \text{in } \Omega \\
\zetabold = \mathbf{0}             & \text{on } \partial \Omega.
\end{cases}
\end{equation}
Due to~\eqref{iso-Delta-Banach}, $\zetabold$
is equal to~$-\Lcalp^{-1} \nablabold q$ in $\Vbfpperp$.
Given $\psibold$ in $\Wbf^{1,\p}_0(\Omega)$,
consider its orthogonal decomposition
$\psibold=\psiboldperp+\psiboldz$, $\psiboldperp$ in~$\Vbfpperp$
and~$\psiboldz$ in~$\Vbfpz$, as in \eqref{orthogonal-Banach}.
For $q$ in $\L^{\pprime}_0(\Omega)$, observe that 
\small\[
\begin{split}
\sup_{\psibold\in\Wbf^{1,\p}_0(\Omega)}  
        \frac{_{\Wbf^{-1,\pprime}}\langle \nablabold \q, \psibold \rangle_{\Wbf^{1,\p}_0}}
        {\SemiNorm{\psibold}_{\Wbf^{1,\p}(\Omega)}}
&   =   \sup_{\psibold\in\Wbf^{1,\p}_0(\Omega)}  
        \frac{_{\Wbf^{-1,\pprime}}\langle\nablabold \q, \psiboldperp \rangle_{\Wbf^{1,\p}_0}}
        {\SemiNorm{\psibold}_{\Wbf^{1,\p}(\Omega)}}
    \le \sup_{\psibold\in\Wbf^{1,\p}_0(\Omega)}  
        \frac{_{\Wbf^{-1,\pprime}}\langle \nablabold \q, \psiboldperp \rangle_{\Wbf^{1,\p}_0}}
        {\SemiNorm{\psiboldperp}_{\Wbf^{1,\p}(\Omega)}}\\
&   =   \sup_{\psiboldperp\in\Vbfpperp} 
        \frac{_{\Wbf^{-1,\pprime}}\langle \nablabold \q, \psiboldperp \rangle_{\Wbf^{1,\p}_0}}
        {\SemiNorm{\psiboldperp}_{\Wbf^{1,\p}(\Omega)}}
\le \sup_{\psibold\in\Wbf^{1,\p}_0(\Omega)}  
        \frac{_{\Wbf^{-1,\pprime}}\langle \nablabold \q, \psibold \rangle_{\Wbf^{1,\p}_0}}
        {\SemiNorm{\psibold}_{\Wbf^{1,\p}(\Omega)}},
\end{split}
\]\normalsize
i.e.,
\begin{equation} \label{equivalence-negative-norms-Banach}
\Norm{\nablabold\q}_{(\Vbfpperp)'}
= \Norm{\nablabold \q}_{\Wbf^{-1,\pprime}(\Omega)}.
\end{equation}
Next, a duality argument reveals
\begin{equation} \label{technical-part-appendix-1-Banach}
\begin{split}
    \Norm{\nablabold q}_{\Wbf^{-1,\pprime}(\Omega)}
&   \overset{\eqref{equivalence-negative-norms-Banach}}{=} 
    \sup_{\Psibold \in\Vbfpperp}
    \frac{_{\Wbf^{-1,\pprime}}\langle\nablabold q, \Psibold\rangle_{\Wbf^{1,\p}_0}}
    {\SemiNorm{\Psibold}_{\Wbf^{1,\p}(\Omega)}}
    \overset{\eqref{elliptic-auxiliary-Banach}}{=}
    \sup_{\Psibold \in\Vbfpperp}
    \frac{_{\Wbf^{-1,\pprime}}\langle - \Lcalp\zetabold, \Psibold\rangle_{\Wbf^{1,\p}_0}}
    {\SemiNorm{\Psibold}_{\Wbf^{1,\p}(\Omega)}} \\
&   \overset{\eqref{Lcalp}}{=}
    \sup_{\Psibold \in\Vbfpperp}
    \frac{_{\Lbb^{\pprime}}\langle\Gcalp\zetabold, \nablau\Psibold\rangle_{\Lbb^{\p}}}
    {\SemiNorm{\Psibold}_{\Wbf^{1,\p}(\Omega)}}
    = \Norm{\Gcalp\zetabold}_{\Lbb^{\pprime}(\Omega)}
    \overset{\eqref{Gcalp}}{=}
    \Norm{\nablau\zetabold}_{\Lbb^{\p}(\Omega)}^{\p-1}\\
&   \overset{\eqref{Lcalp}}{=}
    (_{\Wbf^{-1,\pprime}}\langle -\Lcalp\zetabold,\zetabold\rangle_{\Wbf^{1,\p}_0})^\frac{\p-1}\p
    \overset{\eqref{elliptic-auxiliary-Banach}}{=} 
    (_{\Wbf^{-1,\pprime}}\langle \nablabold q, -\Lcalp^{-1} \nablabold q \rangle_{\Wbf^{1,\p}_0})^\frac{\p-1}\p 
    \\
    & = (_{\L^\p}\langle \div \Lcalp^{-1} \nablabold q , q \rangle_{\L^{\pprime}})^\frac{1}{\pprime}
    \overset{\eqref{Cosserat-operators-Banach}}
    = (_{\L^\p}\langle \vert \Scalstar q\vert^{\pprime-2}\Scalstar q, q\rangle_{\L^{\pprime}})^\frac{1}{\pprime}.
\end{split}
\end{equation}\normalsize

\paragraph*{Showing $\sigma(\Scalpstar,\p,\Gamma,\Omega) = \CNL{\p}{\Gamma}{\Omega}^{-\pprime}$.}
Considering the nonlinear Rayleigh quotient of \eqref{cosserat.eigen2}, we deduce
\small\begin{equation} \label{sigmaSstar=CNL-pprime}
\sigma(\Scalpstar,\p,\Gamma,\Omega)
= \inf_{q\in L^{\pprime}_0(\Omega)}
    \frac{_{\L^{\p}}\langle\vert \Scalstar q\vert^{\pprime-2}\Scalstar q , q \rangle_{\L^{\pprime}}}{\Norm{q}_{\L^{\pprime}(\Omega)}^{\pprime}}
= \inf_{q\in L^{\pprime}_0(\Omega)} \frac{\Norm{\nablabold q}_{\Wbf^{-1,\pprime}(\Omega)}^{\pprime}}{\Norm{q}_{\L^{\pprime}(\Omega)}^{\pprime}}
=: \CNL{\p}{\Gamma}{\Omega}^{-\pprime}.
\end{equation}\normalsize

\paragraph*{Summarizing.}
The assertion~\eqref{homogeneous-BCs-identities-Banach}
follows by noting that
\[
\begin{aligned}
\betaz{\p}{\Gamma}{\Omega}^{\pprime} 
&   \overset{\eqref{beta=CNLmo}}{=}
    \CNL{\p}{\Gamma}{\Omega}^{-\pprime}
    \overset{\eqref{sigmaSstar=CNL-pprime}}{=}
    \sigma(\Scalpstar,\p,\Gamma,\Omega)
    \overset{\eqref{sigmaSp=sigmaSpstar}}{=}
    \sigma(\Scalp,\p,\Gamma,\Omega)^{\pprime-1} \\
&   \overset{\eqref{indices}}{=}
    \sigma(\Scalp,\p,\Gamma,\Omega)^\frac1{\p-1}
    \overset{\eqref{sigmaSp=CBA-p}}{=}
    \CBA{\p}{\Gamma}{\Omega}^{-\frac\p{\p-1}}
    \overset{\eqref{indices}}{=}
    \CBA{\p}{\Gamma}{\Omega}^{-\pprime}.
\end{aligned}
\]

\section{$\W^{1,\p}$-stable right-inverse of the divergence on domains admitting a simplicial tessellation}
\label{appendix:RI-simplices}
Throughout the paper, we assumed that the domain~$\Omega$ is
star-shaped with respect to a ball, cf.~\eqref{rho}.
The only place where this assumption was used
is the proof of Lemma~\ref{lemma:Galdi}
and Theorem~\ref{theorem:new-BA},
where that assumption was required
in order to get explicit bounds on the constants;
cf. \eqref{Galdi's-bound} and~\eqref{constants:Babu-Aziz-mixed}.
In fact, giving up the explicit knowledge of the constant
appearing in the bound,
one can prove these results on more general domains satisfying
\begin{equation} \label{more-general-domains}
\text{$\Omega$ admits a $\gammatilde$ shape-regular decomposition
into~$\Nfraktilde$ simplices.}
\end{equation}
Domains with nontrivial topology
(such as domains with holes, tunnels, \dots)
can satisfy assumption~\eqref{more-general-domains}.
One has the following result.
\begin{theorem}[$\W^{1,\p}$-stable right-inverse of the divergence:
domains admitting a simplicial tessellation]
\label{theorem:Babu-Aziz-general-domain}
Let~$\p$ be in~$(1,\infty)$,
and $\Omega$ and the corresponding
simplicial mesh in~\eqref{more-general-domains}
as above.
There exists a positive constant~$\widetilde C_{BA}(\p,\gammatilde)$
such that for all $f$ in $\L^p_0(\Omega)$,
it is possible to construct~$\vbf$ in~$\Wbf^{1,p}_{0}(\Omega)$
satisfying
\begin{equation} \label{BA-fully-homogeneous-BCS:general-domains}
\div\vbf=f ,
\qquad\qquad\qquad
\SemiNorm{\vbf}_{\Wbf^{1,\p}(\Omega)}
\le \widetilde C_{BA}(\p,\gammatilde, \Nfraktilde) \Norm{f}_{\L^\p(\Omega)}.
\end{equation}
Furthermore,
there exists a positive constant~$\widetilde C_{BA}(\p,\GammaN,\Omega,\gammatilde)$
such that, for all $f$ in $\L^\p(\Omega)$,
it is possible to construct~$\vbf$ in~$\Wbf^{1,\p}_{\GammaN}(\Omega)$
satisfying
\begin{equation} \label{BA-partial-homogeneous-BCS:general-domains}
\div\vbf=f,
\qquad\qquad\qquad
\SemiNorm{\vbf}_{\Wbf^{1,\p}(\Omega)}
\le \widetilde C_{BA}(\p,\GammaN,\Omega,\gammatilde,\Nfraktilde)
        \Norm{f}_{\L^\p(\Omega)}.
\end{equation}
\end{theorem}
\begin{proof}
\noindent \textbf{Proof of~\eqref{BA-fully-homogeneous-BCS:general-domains}}.
Consider the partition of unity~$\{ \phinu \}$,
i.e., for each vertex~$\nu$
of the simplicial mesh
in~\eqref{more-general-domains},
the function~$\phinu$ is the lowest order Lagrangian
function associated with~$\nu$;
the corresponding patch of elements,
i.e., the support of~$\phinu$,
is denoted by~$\omeganu$.
The $\gammatilde$ shape-regularity implies that
there exists a maximum number of simplices
that may be contained in each~$\omeganu$.
As such, there exists a given set of reference patches
$\{\omegatildej\}_{j=1}^{N_{\gammatilde}}$
such that each patch~$\omeganu$
can be mapped via piecewise affine,
global $\W^{1,\infty}$ bijections
into an~$\omegatildej$,
which is star-shaped with respect to a ball.

For each vertex~$\nu$, we construct~$\vbfnu$
in~$\Wbf^{1,\p}_0(\omeganu)$
as in Lemma~\ref{lemma:Galdi}, i.e.,
\[
\div\vbfnu 
= f\phinu -\frac{1}{\vert\omeganu\vert} \int_{\omeganu} f \phinu ,
\]
with stability estimate
\begin{equation} \label{stability-patches}
\SemiNorm{\vbfnu}_{\Wbf^{1,\p}(\omeganu)}
\le C_{\nu}(\p,\gammatilde) \Norm{f}_{\L^\p(\omeganu)}.
\end{equation}
Define $\vbf := \sum_{\nu} \vbfnu$,
where~$\vbfnu$ are extended by zero outside~$\omeganu$.
We have $\vbf_{|\partial\Omega} = \mathbf 0$,
\[
\SemiNorm{\vbf}_{\Wbf^{1,\p}(\Omega)}
\le \sum_{\nu} \SemiNorm{\vbf}_{\Wbf^{1,\p}(\omeganu)}
\overset{\eqref{stability-patches}}{\le}
    (\max_{\nu} C_\nu(\p,\gammatilde))
    \sum_\nu \Norm{f}_{\L^\p(\omeganu)}
\le C_\aleph(\p,\gammatilde, \Nfraktilde)
    \Norm{f}_{\L^\p(\Omega)},
\]
and, due to the properties of the partition of
unity and the fact that $f$ has zero average
over~$\Omega$,
\[
\div \vbf
= (\sum_{\nu} f\phinu)
    - \sum_\nu
    \big(\frac{1}{\vert \omeganu \vert}
        \int_{\omeganu} f\phinu    \big)
    = f - \underbrace{\sum_\nu \big(\frac{1}{\vert \omeganu \vert}
        \int_{\omeganu} f\phinu    \big)}_{\in\Rbb} .
\]
Using that $\vbf$ has zero trace over~$\partial\Omega$
and $f$ has zero average over~$\Omega$ implies
that the last term on the right-hand side is zero,
i.e., $\vbf$ is a $\W^{1,\p}$-stable right-inverse of the divergence of~$f$.

\medskip
\noindent \textbf{Proof of~\eqref{BA-partial-homogeneous-BCS:general-domains}.}
It suffices to use~\eqref{BA-fully-homogeneous-BCS:general-domains}
and arguments as those in Theorem~\ref{theorem:new-BA},
which give the extra dependence of the constant
on~$\GammaN$ and~$\Omega$.
\end{proof}

\begin{remark} \label{remark:comparison-Duran}
The dependence on~$\Nfraktilde$ of the constants
in~\eqref{BA-fully-homogeneous-BCS:general-domains}
and~\eqref{BA-partial-homogeneous-BCS:general-domains}
cannot be removed; cf. the counterexample in \cite[Sect. 3]{Duran:2012}.
\eremk
\end{remark}

\end{document}